\documentclass{article}

\usepackage{a4wide,amsthm, amssymb}

\newtheorem{theorem}{Theorem}[section]
\newtheorem{lemma}[theorem]{Lemma}
\newtheorem{proposition}[theorem]{Proposition}
\newtheorem{corollary}[theorem]{Corollary}
\newtheorem{remark}[theorem]{Remark}
\newtheorem{definition}[theorem]{Definition}

\title{On the cycling operation in braid groups}
\author{Juan Gonz\'{a}lez-Meneses\footnote{Partially supported by
MTM2004-07203-C02-01 and FEDER.} \and Volker Gebhardt}
\date{March 19, 2007}

\begin{document}

\maketitle

\begin{abstract}
The cycling operation is a special kind of conjugation that can be
applied to elements in Artin's braid groups, in order to reduce
their length. It is a key ingredient of the usual solutions to the
conjugacy problem in braid groups. In their seminal paper on
braid-cryptography, Ko, Lee et al. proposed the {\it cycling
problem} as a hard problem in braid groups that could be interesting
for cryptography.  In this paper we give a polynomial solution to
that problem, mainly by showing that cycling is surjective, and
using a result by Maffre which shows that pre-images under cycling
can be computed fast. This result also holds in every Artin-Tits
group of spherical type.

On the other hand, the conjugacy search problem in braid groups is
usually solved by computing some finite sets called  (left) {\it
ultra summit sets} (left-USS), using left normal forms of braids.
But one can equally use right normal forms and compute right-USS's.
Hard instances of the conjugacy search problem correspond to
elements having big (left and right) USS's. One may think that even
if some element has a big left-USS, it could possibly have a small
right-USS.  We show that this is not the case in the important
particular case of {\it rigid} braids. More precisely, we show that
the left-USS and the right-USS of a given rigid braid determine isomorphic
graphs, with the arrows reversed, the isomorphism being defined
using iterated cycling. We conjecture that the same is true for every
element, not necessarily rigid, in braid groups and Artin-Tits
groups of spherical type.
\end{abstract}

\section{Introduction}

Braid groups~\cite{Artin} were related to cryptography in two
independent seminal papers~\cite{AAG,KL}. In both papers, the
security of the proposed cryptosystems relied on the presumed
difficulty of some problems in non-commutative groups, namely the
conjugacy search problem (CSP) and the multiple simultaneous
conjugacy problem (MSCP). They proposed Artin braid groups as good
candidates to implement their cryptosystem, and a lot of literature
has been produced on this subject since then. The results in this
paper refer to braid groups as the main example, but some of them also hold
in other instances of the so-called {\it Garside} groups~\cite{DP,Dehornoy},
which is a family of groups sharing some basic algebraic properties
with braid groups, and which contain all Artin-Tits groups of spherical type.

It seems clear that the main objection to the above cryptosystems,
either in braid groups or in other groups, is the choice of keys. If
one just chooses public and secret keys at random in a braid group,
with given parameters such as length or number of strands, none of
the above cryptosystems can be considered to be secure. It is then
crucial to be able to choose hard instances that resist all known
attacks.

% --
%There are other presumably hard problems in braid groups that have
%been proposed as being possibly interesting for cryptography.
%In~\cite{KL} it is proposed, among others, the {\it cycling
%problem}.  It can be explained as follows. In braid groups one has a
There are other presumably hard problems in braid groups that have
been proposed as being possibly interesting for cryptography.
In~\cite{KL}, the {\it cycling problem}, among others, was suggested.
It can be explained as follows. In braid groups one has a
% -- VG
well known {\it left normal form}, that is, a unique way to write a
braid on $n$ strands $x\in B_n$ as a product $x=\Delta^p x_1\cdots
x_r$, where $\Delta$ is the Garside element, and each $x_i$ is a
simple braid. This normal form will be explicitly defined later. If
we define the {\bf initial factor} of $x$ as $\iota(x)= \Delta^{p}
x_1 \Delta^{-p}$ for $r>0$, and $\iota(x)=1$ for $r=0$, then one has
$x= \iota(x)\: \Delta^p x_2\cdots x_r$. The {\bf left cycling} of
$x$ is defined to be the conjugate of $x$ by its initial factor.
That is, $\mathbf c_L(x)=\Delta^p x_2\cdots x_r \: \iota(x)$.  The
same definition makes sense in every Garside group.

The {\bf cycling problem} asks for, given a braid $y$ and a positive
integer $t$ such that $y$ is in the image of $\mathbf c_L^t$, find a
braid $x$ such that $\mathbf c_L^t(x)=y$.

In this paper we will show that the cycling problem has a polynomial
solution. Namely, it was shown in~\cite{Maffre} that the cycling
problem for $t=1$ has a very efficient solution. That is, if $y$ is
the cycling of some braid, then one can find $x$ such that $\mathbf
c_L(x)=y$ very fast.  In the first part of this paper we will show
the following result, which holds in a special kind of Garside
groups (for instance, it holds in every braid group, and in every
Artin-Tits group of spherical type).

\begin{theorem}\label{T:surjective}
If $G$ is a Garside group which is atom-friendly (on the left), then
$\mathbf c_L : \: G \rightarrow G$ is surjective.
\end{theorem}

As an immediate corollary, a solution to the cycling problem is just
given by applying $t$ times the algorithm in~\cite{Maffre}. This
clearly gives a polynomial solution to the cycling problem, since it
is so for $t=1$.

The proof of Theorem~\ref{T:surjective} makes use not only of left
normal forms, but of {\it right normal forms} of elements in $B_n$
(or in $G$).  We shall see that, under certain conditions, an
inverse of $x$ under cycling, using left normal forms, is precisely
the cycling of $x$ using right normal forms. This shows that {\it
left} and {\it right cyclings}, $\mathbf c_L$ and $\mathbf c_R$, are
closely related.

The cycling operation is mainly used to find simpler conjugates of a
braid, and also to compute finite sets which are invariants of
conjugacy classes and allow to solve the conjugacy problem in $B_n$.
One of such sets is the ultra summit set of a given braid $x$,
$USS(x)$.  One usually defines this set by using left normal forms,
but it is equally possible to define it using right normal forms,
hence one usually has two finite sets associated to $x$, that we
denote $USS_L(x)$ and $USS_R(x)$.

% --
% The algorithm by Gebhardt~\cite{Gebhardt} to solve the conjugacy
%search problem in braid groups (and in any Garside group) relies on
The algorithmic solution to the conjugacy search problem in braid
groups (and in any Garside group) developed in \cite{Gebhardt} relies on
% -- VG
computing ultra summit sets. Hence braids having small ultra summit
sets are not hard instances for the conjugacy search problem. This
means that if one wants to find a good key for a cryptographic
protocol, one needs to choose a braid with a big ultra summit set.
But we have seen that there are two kind of ultra summit sets,
$USS_L(x)$ and $USS_R(x)$, and the question arises on whether one of
them can be big while the other one is small.

On the other hand, there are three geometric kind of braids:
periodic, reducible and pseudo-Anosov~\cite{BLM}. The conjugacy
search problem for periodic braids is solvable in polynomial
time~\cite{BGG3}. Reducible braids are those which can be
decomposed, in some sense, into braids with fewer strands. There are
algorithms to find this decomposition~\cite{BGN}, see
also~\cite{Lee_reducible}, although they are not polynomial.
Nevertheless, in most cases the decomposition can be found very
fast, and the conjugacy problem is split into several conjugacy
problems on fewer strands.  Hence, it would be desirable to know
pseudo-Anosov braids whose ultra summit sets are big.

But one can solve the conjugacy search problem for pseudo-Anosov
braids using {\it rigid} braids (these will be defined later):
In~\cite{GM_roots} it is shown that the conjugacy search problem for
two pseudo-Anosov braids $x$ and $y$ is equivalent to the same
problem for $x^m$ and $y^m$, for every nonzero integer $m$. And
in~\cite{BGG1} it is shown that every pseudo-Anosov element in its
ultra summit set, has a small power which is rigid (we will be more
explicit in the next section).  Therefore, one just needs to care
about rigid braids.  So the above question is transformed into the
following: if $x$ is a rigid braid, is it possible that $USS_L(x)$
is big and $USS_R(x)$ is small, or vice versa?   The answer is
negative, and it is given by the following results.

\begin{theorem}\label{T:rigid_USS}
A braid $x\in B_n$ with $\ell(x)>1$ is conjugate to a left rigid braid if and only if it is conjugate to a right rigid braid.
\end{theorem}

In the above case, we will show that $\#(USS_L(x))=\#(USS_R(x))$.
Therefore, if one is able to find a rigid element $x$ such that
$USS_L(x)$ is big, the same will happen with $USS_R(x)$, so the
conjugacy search problem will be equally difficult by using either
left or right normal forms.

Moreover, we will show that the relation between $USS_L(x)$ and $USS_R(x)$ is deeper than just having the same number of elements. In order to compute $USS_L(x)$ using the algorithm in~\cite{Gebhardt}, one actually computes a directed graph, that we will denote $USG_L(x)$ (left ultra summit graph of $x$). The vertices of $USG_L(x)$ correspond to the elements of $USS_L(x)$, and the arrows are labeled by simple braids, in such a way that there is an arrow labeled by $s$, going from $u$ to $v$, if and only if $s^{-1}us=v$.  In the same way, one can define $USG_R(x)$, where in this case the vertices correspond to elements in $USS_R(x)$, and there is an arrow labeled by $s$, going from $u$ to $v$, if and only if $ sus^{-1}=v$.  We will denote by $USG_R(x)^{op}$ the graph which is isomorphic to $USG_R(x)$ as a (non-directed) graph, but with the arrows reversed. The result that compares the graphs $USS_L(x)$ and $USS_R(x)$ is the following:

\begin{theorem}\label{T:rigid_USG}
Let $x\in B_n$ with $\ell(x)>1$ be conjugate to a left rigid braid. Then $USG_L(x)$ and $USG_R(x)^{op}$ are isomorphic directed graphs.
\end{theorem}

\begin{remark}
We recently learnt from Jean Michel, Fran\c{c}ois Digne et David Bessis, that $USG_L(x)$ (and thus $USG_R(x)$) are Garside categories. In this context, the notation $USG_R(x)^{op}$ makes sense, since it refers to the opposite category. Then Theorem~\ref{T:rigid_USG} says that $USG_L(x)$ and $USG_R(x)^{op}$ are isomorphic Garside categories. Or in other words, there exists a contravariant isomorphism from $USG_L(x)$ to $USG_R(x)$
\end{remark}

This paper is structured as follows: In Section~\ref{S:Garside} some basic notions of braids and Garside theory are given. Specialists in Garside theory may skip this Section and go directly to Section~\ref{S:surjective}, in which Theorem~\ref{T:surjective} is shown.  The proofs of Theorems~\ref{T:rigid_USS} and \ref{T:rigid_USG} are given in Section~\ref{S:rigid_USS}.

{\bf Acknowledgements:}  This paper was conceived in the framework of a collaboration of the authors with Joan S. Birman. Most arguments in it have been discussed with her, and in particular she participated in finding the right conjectures that became Theorems~\ref{T:rigid_USS} and \ref{T:rigid_USG}.  We are deeply grateful to her for these contributions, and also for her advice and support.  The first author thanks Thierry Berger and Samuel Maffre for inviting him to Limoges, to the PhD defense of the latter, where he learnt about the results which are a key tool in Section~\ref{S:surjective}.

\section{Basic ingredients of Garside theory.}\label{S:Garside}

In this section we will explain the notions and results that will be
used throughout the rest of the paper. Namely, we will briefly
describe the basic ingredients of the Garside structure of braid
groups. In general, a Garside group is a group satisfying the
structural properties defined in this section, and the main examples
are braid groups and Artin-Tits groups of spherical type. For a
short introduction to Garside theory, with a precise definition of a
Garside group, see~\cite{BGG1}.

The braid group on $n$ strands, $B_n$ can be defined by its well
known group presentation~\cite{Artin}:
$$
   B_n= \left\langle  \sigma_1,\ldots, \sigma_{n-1} \left|
   \begin{array}{ll} \sigma_i\sigma_j=\sigma_j\sigma_i, & \mbox{ if } |j-i|>1 \\
    \sigma_i\sigma_j\sigma_i=\sigma_j\sigma_i\sigma_j, & \mbox{ if } |j-i|=1
   \end{array} \right. \right\rangle.
$$
If we consider the above as a monoid presentation, this defines the
monoid of positive braids, $B_n^+$.  Garside~\cite{Garside} showed
that $B_n^+$ embeds into $B_n$, so the elements of $B_n^+$, called
{\it positive braids} can be seen as the braids in $B_n$ that can be
written as a word in the generators (but not their inverses).  There
is a special positive element, called {\it half twist} or {\bf
Garside element}, defined by $\Delta = \sigma_1 (\sigma_2 \sigma_1)
\cdots (\sigma_{n-1}\cdots \sigma_1)$. Artin~\cite{Artin} showed
that the center of $B_n$ is the cyclic subgroup generated by
$\Delta^2$. In general, every Garside group has a distinguished
monoid of positive elements, and a special Garside element,
% --
%$\Delta$, which has a central power.
$\Delta$, which has a central power $\Delta^e$.
Conjugation by $\Delta$ is an inner automorphism which preserves the
set of simple elements; we denote this automorphism by $\tau$.
% -- VG

In $B_n$ one can define two partial relations, related to left and
right divisibility, respectively. Namely, given $a,b\in B_n$ we say
that $a\preccurlyeq b$ if $a^{-1}b\in B_n^+$, that is, if $ap=b$ for
some positive braid $p$. We then say that $a$ is a {\it
left-divisor}, or a {\it prefix} of $b$. On the other hand, we say
that $a\succcurlyeq b$ if $ab^{-1}\in B_n^+$, that is, if $a=pb$ for
some positive braid $p$. In this case we say that $b$ is a {\it
right-divisor}, or a {\it suffix} of $a$. Notice that $B_n^+ =
\{p\in B_n; \; 1\preccurlyeq p\} = \{p\in B_n;\; p\succcurlyeq 1\}$.

Each of the above partial orders define a lattice structure on
$B_n$. This means that given two braids $a,b\in B_n$, there exist a
unique greatest common divisor $a \wedge_L b $ and a unique least
common multiple $a\vee_L b$, naturally defined by the left
divisibility relation $\preccurlyeq$, and also unique gcd's and
lcm's, $a\wedge_R b$ and $a\vee_R b$, naturally defined by
$\succcurlyeq$.

In $B_n$, the generators $\sigma_1,\cdots,\sigma_{n-1}$ are called
{\bf atoms}. In general, in a Garside group, an atom is a positive
element that cannot be decomposed as a product of two positive
elements. In the particular case of $B_n$ and of Artin-Tits groups
of spherical type, the Garside element $\Delta$ is the (left and
right) least common multiple of all atoms. This is not true in
general for other Garside groups, and this is one of the reasons why
the proof of Theorem~\ref{T:surjective} above does not generalize to
every Garside group.

Several normal forms for elements in $B_n$ have been defined. We
will concentrate in the one defined independently by
Adjan~\cite{Adjan}, Deligne~\cite{Deligne}, Elrifai-Morton~\cite{EM}
and Thurston~\cite{Epstein}, which is an improvement of the solution
to the word problem given by Garside~\cite{Garside}. We say that a
braid is {\bf simple} if it is a positive prefix of $\Delta$. It is
well known that this happens if and only if it is a positive {\it
suffix} of $\Delta$. The set of simple braids is then $S=
\{s \in B_n;\; 1\preccurlyeq s \preccurlyeq \Delta\} = \{s\in B_n;\;
\Delta \succcurlyeq s \succcurlyeq 1\}$.

% --
%Given two simple elements $s,s'$, we say that the decomposition
%$ss'$ is {\bf left-weighted} if $s$ is the biggest simple prefix of
%$ss'$, in the sense that every simple prefix of $ss'$ is a prefix of
%$s$.  This can also be expressed by $s=(ss')\wedge_L \Delta$.  We
%say that $ss'$ is {\bf right-weighted} if $s'$ is the maximal simple
%suffix of $ss'$, in other words, $s'= (ss')\wedge_R \Delta$.
\begin{definition}
\label{D:left_right_weighted}
Given two simple elements $s,s'$, we say that the decomposition
$ss'$ is {\bf left-weighted} if $s$ is the maximal simple prefix of
$ss'$, that is, if $s=(ss')\wedge_L \Delta$.  Similarly, we say that
$ss'$ is {\bf right-weighted} if $s'$ is the maximal simple
suffix of $ss'$, that is, if $s'= (ss')\wedge_R \Delta$.

For a simple element $s$ we call $\partial (s) = s^{-1}\Delta$ the
{\bf right complement} of $s$.  Note that as $s\preccurlyeq\Delta$ and
$s\,\partial (s)=\Delta$, the element $\partial (s)$ is simple.
Hence, this defines a map $\partial: S \rightarrow S$ on the set $S$
of simple elements. As
$\partial(\partial(s))=\Delta^{-1}s\Delta=\tau(s)$ for any simple $s$,
the map $\partial$ is a bijection on $S$ and $\partial^2 = \tau$.
We similarly define the {\bf left complement} of $s$ as $\Delta
s^{-1} = \Delta\partial(s)\Delta^{-1} = \tau^{-1}(\partial(s)) =
\partial^{-1}(s)$.
\end{definition}

Observe that, given two simple elements $s$ and $s'$, the
product $ss'$ is left weighted if and only if there is no prefix
$t\preccurlyeq s'$ such that $st$ is simple, or in other words, such
that $t\preccurlyeq\partial(s)$. Hence $ss'$ is left weighted if and
only if $\partial(s)\wedge_L s' =1$.
Similarly, $ss'$ is right weighted if and only if
$s \wedge_R \partial^{-1}(s') =1$.
% -- VG

\begin{definition}
Given a braid $x\in B_n$, its {\bf left normal form} is a
decomposition $x=\Delta^p x_1\cdots x_r$, satisfying the following
conditions:
\begin{enumerate}

 \item $p\in \mathbb Z$ is the maximal integer such that $\Delta^{-p}x$ is positive.

 \item $x_i= (x_i\cdots x_r)\wedge_L \Delta \neq 1$ for $i=1,\ldots, r$.

\end{enumerate}
\end{definition}

In other words, each $x_i$ is a proper simple element (different
from $1$ and $\Delta$), and it is the biggest simple prefix of
$x_i\cdots x_r$.  It is well known that normal forms can be
recognized `locally'. This means that $\Delta^{p}x_1\cdots x_r$ is
in left normal form if and only if each $x_i$ is a proper simple
element and $x_ix_{i+1}$ is left-weighted for $i=1,\ldots, r-1$.
The left normal form of a braid exists and it is unique.  The
integers $p$ and $r$ are then determined by $x$, so one can define
the {\bf infimum}, {\bf supremum} and {\bf canonical length} of $x$,
respectively, by $\inf(x)=p$, $\; \sup(x)=p+r$ and $\ell(x)=r$.
This terminology is explained by noticing that $p$ and $p+r$ are,
respectively, the biggest and the smallest integers such that
$\Delta^p \preccurlyeq x \preccurlyeq \Delta^{p+r}$, which is
usually written $x\in [\Delta^{p},\Delta^{p+r}]$, or simply $x\in
[p,p+r]$. The canonical length $r$ is just the size of this
interval, which corresponds to the number of non-Delta factors in
the left normal form of $x$.

We notice that one has the analogous definitions related to
$\succcurlyeq$:

\begin{definition}
Given a braid $x\in B_n$, its {\bf right normal form} is a
decomposition $x= y_1\cdots y_r\Delta^p$, satisfying the following
conditions:
\begin{enumerate}

 \item $p\in \mathbb Z$ is the maximal integer such that $x\Delta^{-p}$ is positive.

 \item $y_i= (y_1\cdots y_i)\wedge_R \Delta \neq 1$ for $i=1,\ldots, r$.

\end{enumerate}
\end{definition}

The property of being a right normal form is also a local property
($y_iy_{i+1}$ is right-weighted for every $i$), and this
decomposition also exists and is unique for each braid. We remark
that the integers $p$ and $r$ in this case are exactly the same as
those corresponding to the left normal form. This means that
$\inf(x)=p$ and $\sup(x)=p+r$ are, respectively, the maximal and
minimal integers such that $\Delta^{p+r}\succcurlyeq x \succcurlyeq
\Delta^p$, hence $\inf(x)$, $\sup(x)$ and $\ell(x)$ can be equally
defined using right normal forms instead of left normal forms.

Recall that we defined the initial factor of a braid in the
introduction. Since we are using two distinct structures in $B_n$,
we will define left and right versions of initial and final factors,
as follows.
%--
%First we notice that conjugation by $\Delta$ is an
%inner automorphism of $B_n$ that preserves the set of simple
%elements, and that we denote by $\tau$.
% --VG
Given $x=\Delta^p x_1\cdots
x_r$ in left normal form, we define its {\bf left initial factor} as
$\iota_L(x)= \tau^{-p}(x_1)$, and its {\bf left final factor} by
$\varphi_L(x)=x_r$.  In the same way, if $x=y_1\cdots y_r \Delta^p$
is in right normal form, we define its {\bf right initial factor} by
$\iota_R(x)=\tau^p(y_r)$, and its {\bf right final factor} by
$\varphi_R(x)=y_1$.

There are special maps from the braid group to itself that consist
of conjugating each element by the above initial or final factors.
These operations, called cyclings and decyclings, are key
ingredients in most of the known solutions to the conjugacy problem
in braid groups. The precise definition is as follows.

\begin{definition}
The following maps, from $B_n$ to itself, are defined for each $x\in
B_n$ as follows:
\begin{enumerate}

\item {\bf Left cycling:}   $\mathbf c_L(x)= \iota_L(x)^{-1}\cdot x \cdot\iota_L(x)$.

\item {\bf Left decycling:}   $\mathbf d_L(x)= \varphi_L(x) \cdot x \cdot \varphi_L(x)^{-1}$.

\item {\bf right cycling:}   $\mathbf c_R(x)= \iota_R(x) \cdot x \cdot \iota_R(x)^{-1}$.

\item {\bf right decycling:}   $\mathbf d_R(x)= \varphi_R(x)^{-1}\cdot x \cdot\varphi_R(x)$.

\end{enumerate}
\end{definition}

In other words, if $x=\Delta^p x_1\cdots x_r$ is in left normal form, then
$$
\mathbf c_L(x)=\Delta^p x_2\cdots x_r \tau^{-p}(x_1), \hspace{1cm}
\mathbf d_L(x)=x_r \Delta^p x_1\cdots x_{r-1},
$$
and if $x= y_1\cdots y_r \Delta^p$ is in right normal form, then
$$
\mathbf c_R(x)=\tau^p(y_r) y_1\cdots y_{r-1} \Delta^p, \hspace{1cm}
\mathbf d_R(x)= y_2\cdots y_r \Delta^p y_1.
$$

We notice that there is an involution of the braid group, $rev :\:
B_n \rightarrow B_n$, which sends every braid $x=
\sigma_{i_1}^{e_1}\cdots \sigma_{i_m}^{e_m}$ to its {\bf reverse}
$rev(x)=\overleftarrow{x} =  \sigma_{i_m}^{e_m}\cdots
\sigma_{i_1}^{e_1}$, that is, the same word read backwards.
% --
%This is
Observe that the map $rev$ is well-defined, as the relations of $B_n$
are invariant under $rev$.  The map $rev$ is
% -- VG
an anti-isomorphism, and one can easily check that the left normal
form of $x$ is mapped by $rev$ to the right normal form of
$\overleftarrow x$, and vice versa. Also $\overleftarrow{\iota_R(x)}
= \iota_L (\overleftarrow{x}) $,
$\overleftarrow{\varphi_R(x)}=\varphi_L(\overleftarrow{x})$, and
then  $\overleftarrow{\mathbf c_R(x)} = \mathbf c_L (\overleftarrow
x)$ and $\overleftarrow{\mathbf d_R(x)} = \mathbf d_L
(\overleftarrow x)$.  This means that applying $\mathbf c_R$ and
$\mathbf d_R$ to a braid $x$ corresponds to applying the usual
cycling and decycling operations, $\mathbf c_L$ and $\mathbf d_L$,
to its reverse $\overleftarrow x$.  This implies that all results
which are usually shown using left normal forms, $\mathbf c_L$ and
$\mathbf d_L$, will also hold using right normal forms, $\mathbf
c_R$ and $\mathbf d_R$, by symmetry.

Cyclings and decyclings have been used to define suitable finite
subsets of $B_n$ which allow to solve the conjugacy decision problem
and the conjugacy search problem in braid groups. Namely, the {\bf
super summit set} of an element $x$, denoted $SSS(x)$~\cite{EM} is
defined as follows. If we denote $C(x)$ the conjugacy class of $x$,
then
$$
   SSS(x)=\{y\in C(x);\quad \ell(y) \mbox{ is minimal}\}.
$$
Notice that this set does not depend on which structure of $B_n$
(left or right) we used to define $\ell(y)$. A subset of $SSS(x)$ is
the {\it ultra summit set} of $x$~\cite{Gebhardt}. In this case,
since $USS(x)$  is defined by using cyclings, one needs to
distinguish between the {\bf left ultra summit set} of $x$,
$$
   USS_L(x)=\{y\in SSS(x);\quad  \exists t\geq 1, \; \mathbf c_L^t(y)=y \},
$$
and the {\bf right ultra summit set} of $x$,
$$
   USS_R(x)=\{y\in SSS(x);\quad  \exists t\geq 1, \; \mathbf c_R^t(y)=y \}.
$$

Both $SSS(x)$, $USS_L(x)$ and $USS_R(x)$ are, by definition,
invariants of the conjugacy class of $x$. Hence one can determine
whether two braids $x,y\in B_n$ are conjugate by computing, say,
$USS_L(x)$ and $USS_L(y)$ and checking if they are equal. Actually,
it suffices to compute $USS_L(x)$, one element $y'\in USS_L(y)$ and
to check whether $y'\in USS_L(x)$. In~\cite{EM} it is shown how to
compute $SSS(x)$, and~\cite{Gebhardt} gives an algorithm to compute
$USS_L(x)$ (which can also be used to compute $USS_R(x)$). More precisely, the algorithm computes a directed graph whose set of vertices is $USS_L(x)$. We will define such a graph as follows.

\begin{definition}
Given $x\in B_n$, we define the {\bf left ultra summit graph} of $x$, denoted $USG_L(x)$, as the directed graph whose set of vertices is $USS_L(x)$ and whose arrows are labeled by simple elements, in such a way that there is an arrow labeled $s$, starting at $u$ and ending at $v$, if $s^{-1} u s = v$.

In the same way, we define the {\bf right ultra summit graph} of $x$, denoted $USG_R(x)$, as the directed graph whose set of vertices is $USS_R(x)$ and whose arrows are labeled by simple elements, in such a way that there is an arrow labeled $s$, starting at $u$ and ending at $v$, if $s u s^{-1} = v$.
\end{definition}

We remark that in~\cite{Gebhardt}, the graph that is computed is not precisely $USG_L(x)$, but one with less arrows:

\begin{definition}
Given $x\in B_n$, we define the graph $minUSG_L(x)$ to be the subgraph of $USG_L(x)$ with the same set of vertices, but only with {\bf minimal} arrows. An arrow labeled by $s$ and starting at $u$  is said to be minimal if it cannot be decomposed as a product of arrows, that is, if there is no directed path in $USG_L(x)$ starting at $u$, with labels $s_1,\ldots, s_k$, such that $s=s_1\cdots s_k$.

In the same way, we define the graph $minUSG_R(x)$ to be the subgraph of $USG_R(x)$ with the same set of vertices, but only with minimal arrows.
\end{definition}

It is known that all the above graphs are connected. The arrows in these graphs allow to know how to connect, by a conjugation, $x$ to any element in $USS_L(x)$ and
$y$ to any element in $USS_L(y)$. Hence, the above procedure also solves
the conjugacy search problem in $B_n$ (and in any Garside group),
that is, it finds a conjugating element from $x$ to $y$ provided it
exists.

In~\cite{BGG1} one can find is a project to find a polynomial solution to
the conjugacy search problem in braid groups. One of the crucial
open problems in this project concerns {\it rigid} braids, which are
defined as follows.  As above, since we are using two different
structures of $B_n$ we will define rigid elements on the left and on
the right. In this way, we will say that an element $x=\Delta^p
x_1\cdots x_r$ (written here in left normal form, with $r>0$) is {\bf left
rigid}, if $\Delta^p x_1 \cdots x_r \iota_L(x)$ is in left normal
form as written.   In the same way, we will say that $x=y_1\cdots
y_r \Delta^p $ (written in right normal form, with $r>0$) is {\bf right rigid}
if $\iota_R(x) y_1\cdots y_r \Delta^p$ is in right normal form as
written, or alternatively, if $\overleftarrow x$ is left rigid.
These are the elements that have the best possible behavior with
respect to cyclings and decyclings, since in this case iterated
cyclings or decyclings just correspond to cyclic permutation of the
factors (for non-rigid elements this is not the case, since one
needs to compute the left normal form of $\mathbf c_L(x)$ in order
to be able to apply $\mathbf c_L$ again, and this modifies some of
the original factors of $x$).

There are some interesting results concerning rigid braids.

\begin{theorem}{\rm \cite{BGG1}}
If $x\in B_n$ is left $[$right$]$ rigid then $x\in USS_L(x)$
$[x\in USS_R(x)]$. Moreover, if $\ell(x)>1$ then  $USS_L(x)$ $[USS_R(x)]$ is
precisely the set of left $[$right$]$ rigid conjugates of $x$.
\end{theorem}

\begin{theorem}{\rm \cite{BGG1}}
If $x\in B_n$ is a pseudo-Anosov braid, and $x\in USS_L(x)$
$[x\in USS_R(x)]$, then $x^m$ is left $[$right$]$ rigid for
some $m< (\frac{n(n-1)}{2})^3$.
\end{theorem}

Since pseudo-Anosov braids seem to be generic in $B_n$, and the
conjugacy search problem for pseudo-Anosov braids $x$ and $y$ can be
solved just by solving it for $x^m$ and $y^m$ for any $m\neq
0$~\cite{GM_roots}, the rigid case turns out to be probably the most
important case to solve the conjugacy search problem in $B_n$.

% --
%As was noticed by Gebhardt in~\cite{Gebhardt}, if the canonical
As was noticed in~\cite{Gebhardt}, if the canonical
% -- VG
length of a random braid $x$ is big enough with respect to the
number of strands, then $USS_L(x)$ consists exactly of $2\ell(x)$
elements in 100\% of the tested cases, meaning that the probability
of getting a larger $USS_L(x)$ seems to tend to zero very rapidly as
$\ell(x)$ grows. Moreover, in this `generic' cases the braids in
$USS_L(x)$ are pseudo-Anosov and left rigid. We remark that
Gebhardt's algorithm is a deterministic algorithm that is
`generically' polynomial, although there is no written proof, to our
knowledge, that either pseudo-Anosov braids or  braids conjugate to
a rigid element are generic in $B_n$.

There are instances of left rigid elements whose ultra summit set is
much bigger than expected. For instance, as is noticed
in~\cite{BGG1}, the braid in $B_{12}$
\begin{eqnarray*}
E &=&(\sigma_2 \sigma_1 \sigma_7 \sigma_6 \sigma_5 \sigma_4 \sigma_3 \sigma_8 \sigma_7 \sigma_{11} \sigma_{10}) \cdot (\sigma_1 \sigma_2 \sigma_3 \sigma_2 \sigma_1 \sigma_4 \sigma_3 \sigma_{10}) \cdot  \\
& & (\sigma_1 \sigma_3 \sigma_4 \sigma_{10}) \cdot  (\sigma_1 \sigma_{10}) \cdot  (\sigma_1 \sigma_{10} \sigma_9 \sigma_8 \sigma_7 \sigma_{11}) \cdot  (\sigma_1 \sigma_2 \sigma_7 \sigma_{11})
\end{eqnarray*}
is a pseudo-Anosov, rigid braid with $\ell(E)=6$, such that
$\#(USS_L(E))=264 = 44\cdot 6$, instead of the expected value of
$12=2\cdot 6$. Also, the braid in $B_{12}$
\begin{eqnarray*}
F & = & (\sigma_3 \sigma_2 \sigma_1 \sigma_4 \sigma_6 \sigma_8 \sigma_7 \sigma_6 \sigma_9 \sigma_{10}\sigma_{11} \sigma_{10}) \cdot (\sigma_1 \sigma_2 \sigma_4 \sigma_3 \sigma_2 \sigma_1 \sigma_5 \sigma_7 \sigma_{10} \sigma_{11} \sigma_{10}) \cdot  \\
& &  (\sigma_3 \sigma_5 \sigma_7
\sigma_{10 }\sigma_{11} \sigma_{10}) \cdot (\sigma_3 \sigma_5 \sigma_7 \sigma_6 \sigma_8 \sigma_{10} \sigma_{11})
\end{eqnarray*}
is pseudo-Anosov and rigid, with $\ell(F)=4$ and $\#(USS_L(F))=232 =
58\cdot 4$, instead of the expected value of $8=2\cdot 4$. The
reason why these special examples of rigid braids exist, and how one
can compute them, is still a mystery. Solving this problem would be
an important step towards finding secure keys for cryptographic
protocols with braid groups.

But recall that we are considering two distinct structures in $B_n$.
Hence it could be possible, a priori, that $USS_R(E)$ or $USS_R(F)$
are much smaller that $USS_L(E)$ or $USS_L(F)$, respectively.
Theorem~\ref{T:rigid_USS} tells us that this is not the case, since
$\#(USS(x))=\#(USS(x))$ for every rigid braid $x$ of canonical
length greater than 1.

\section{Cycling is surjective}\label{S:surjective}

In this section we will show Theorem~\ref{T:surjective}, that is, we
will show that $\mathbf c_L$ (and thus $\mathbf c_R$) is a
surjective map.

% --
%First we recall that, for every simple element $s$, by definition
%one has $s\preccurlyeq \Delta$, meaning that there exist a positive
%element $s^*$ such that $s\: s^*=\Delta$. Since $s^*$ is then a
%positive right divisor of $\Delta$, it follows that $s^*$ is simple.
%We say that $s^*$ is the {\bf right complement} of $s$. This defines
%a map $\partial: S \rightarrow S$, by $\partial (s)=s^*$.  One can
%see that $\partial$ is a bijection of the set of simple elements
%$S$, and that $\partial^2 = \tau$, where $\tau$ is the inner
%automorphism determined by $\Delta$.
%
%We also note that, given two simple elements $s$ and $s'$, the
%product $ss'$ is left weighted if and only if there is no prefix
%$t\preccurlyeq s'$ such that $st$ is simple. In other words, if
%there is no nontrivial prefix of $s'$ and $\partial(s)$. Hence $ss'$
%is left-weighted if and only if $\partial(s)\wedge_L s' =1$.
First we recall the definition of the right complement $\partial (s)$
of a simple element $s$ from Definition~\ref{D:left_right_weighted}.
A product $ss'$ of two simple elements $s$ and $s'$ is left-weighted
if and only if $\partial(s)\wedge_L s' =1$.
% -- VG

It was shown by Maffre~\cite{Maffre} that the pre-image of a braid
$x\in B_n$ under $\mathbf c_L$ can be computed fast, provided that
$x$ is in the image of $\mathbf c_L$. The procedure depends on
whether the infimum of the existing pre-image of $x$ is equal to $\inf(x)$ or not. We will treat the situation from a slightly different point of view, although the pre-images that we will compute are exactly the same as those given by Maffre.

The following result holds for every Garside group $G$. In the
particular case of $B_n$, recall that the atoms are just the
generators $\sigma_1,\ldots,\sigma_{n-1}$. We will see that in some
particular cases, we can obtain a pre-image of $x$ by $\mathbf c_L$,
just by conjugating $x$ by an atom, and then by $\Delta^{-1}$.

\begin{proposition}\label{P:preimage1}
Let $G$ be a Garside group, and let $x=\Delta^p x_1\cdots x_r\in G$
be written in left normal form. If there is an atom $a$ such that
$\tau^p(a)\not\preccurlyeq x_1\cdots x_r a$, then $\mathbf
c_L(\tau^{-1}(a^{-1}xa))=x$.
\end{proposition}

\begin{proof}
Define $z=a^{-1}xa=\partial(a) \Delta^{p-1} x_1\cdots x_r a =
\Delta^{p-1} \partial^{2p-1}(a) x_1\cdots x_r a$. Notice that
$\partial(\partial^{2p-1}(a)) = \partial^{2p}(a) =
\tau^p(a)\not\preccurlyeq x_1\cdots x_r a$.  But $\tau$ transforms
atoms into atoms, hence $\tau^p(a)$ is an atom. This means that
$\tau^p(a)\not\preccurlyeq x_1\cdots x_r a$ is equivalent to
$\tau^p(a)\wedge_L x_1\cdots x_r a=1$, since an atom has no
nontrivial prefixes.

Notice that $\Delta\not\preccurlyeq  x_1\cdots x_ra$, otherwise
$a\preccurlyeq \Delta \preccurlyeq x_1\cdots x_r a$. Hence
$\inf(x_1\cdots x_ra)=0$ which implies that $\iota(x_1\cdots x_ra)$
is precisely the biggest simple prefix of $x_1\cdots x_ra$.
Therefore, since $\tau^p(a)\wedge_L x_1\cdots x_r a=1$, we also have
$\tau^p(a)\wedge_L\iota(x_1\cdots x_r a)=1$. In other words, if
$z_2\cdots z_k$ is the left normal form of $x_1\cdots x_r a$, then
$\tau^p(a)\wedge_L z_2=1$, that is $\partial(\partial^{2p-1}(a))
\wedge_Lz_2=1$, so $\partial^{2p-1}(a) z_2$ is left-weighted.  This
implies that $\partial^{2p-1}(a) z_2\cdots z_k$ is the left normal
form of $\partial^{2p-1}(a) x_1\cdots x_r a$. Hence $\iota(z)=
\tau^{-p+1}(\partial^{2p-1}(a)) =
\partial^{-2p+2}(\partial^{2p-1}(a))= \partial(a)$.

If we apply left-cycling to $z$, we then obtain
$$
\mathbf c_L (z)= z^{\partial(a)}=\Delta^{p-1} x_1\cdots x_r a
\partial(a) = \Delta^{p-1} x_1\cdots x_r \Delta = \tau(x)
$$
It is well known (and can be derived from the definitions and from
the fact that $\tau$ is a bijection of $S$) that $\tau$ sends left
(and right) normal forms to left (and right) normal forms. Hence
$\tau$ commutes with $\mathbf c_L$ (and with $\mathbf c_R$).
Therefore $\mathbf c_L(\tau^{-1}(z))=\tau^{-1}(\mathbf
c_L(z))=\tau^{-1}(\tau(x))=x$, as we wanted to show.
\end{proof}

We will now see that, in the cases where the hypothesis of
Proposition~\ref{P:preimage1} are not satisfied, then a preimage by
$\mathbf c_L$ of $x$ is just $\mathbf c_R(x)$.  This time our proof
does not work for every Garside group, but we need some special
property to be satisfied.  Given a Garside group $G$, we will denote
by $\mathcal A$ the set of atoms. Given a simple element $s\in G$,
we will define the {\bf starting set} of $s$ as $\mathcal
S(s)=\{a\in \mathcal A;\; a\preccurlyeq s\}$.

\begin{definition}
Given a Garside group $G$, we will say that $G$ is {\bf
atom-friendly} (on the left) if
\begin{enumerate}

 \item $\mbox{\rm lcm}_L(\mathcal A)=\Delta$.

 \item $\mathcal S(\mbox{\rm lcm}_L(\mathcal B)) =\mathcal B $ for every $\mathcal B\subset \mathcal A$.

\end{enumerate}
\end{definition}

We remark that the terminology {\it atom-friendly} is new. To our
knowledge, no common name has been given to those Garside groups
satisfying the above two conditions. It is nevertheless well
known~\cite{Picantin} that braid groups, and more generally
Artin-Tits group of spherical type are atom-friendly (on the left
and on the right). Hence the following result holds in all
Artin-Tits groups of spherical type.

\begin{proposition}\label{P:preimage2}
Let $G$ be a Garside group which is atom-friendly (on the left). Let
$x=\Delta^p x_1\cdots x_r\in G$ be written in left normal form. If
for every atom $a$ one has $\tau^p(a)\preccurlyeq x_1\cdots x_r a$,
then $\mathbf c_L(\mathbf c_R(x)))=x$.
\end{proposition}

\begin{proof}
Let us define ${\cal D}$ to be the set of atoms $a$ such that
$\tau^p(a)\not \preccurlyeq x_1$. That is ${\cal D}={\cal
A}\backslash {\cal S}(\tau^{-p}(x_1))= \mathcal A \backslash
\mathcal S(\iota(x))$. Define also the simple element
$D=\mbox{lcm}_L(\mathcal D)$.  Let us show that $\Delta \preccurlyeq
x_1 \cdots x_r D$.  Indeed, for every atom $a\notin {\cal D}$ one
has $\tau^p(a)\preccurlyeq x_1 \preccurlyeq x_1\cdots x_r D$, and
for every atom $a\in {\cal D}$ one has $a\preccurlyeq D$, so using
the hypothesis it follows that $\tau^p(a)\preccurlyeq x_1\cdots x_r
a \preccurlyeq x_1\cdots x_r D$. Therefore $\tau^p (a) \preccurlyeq
x_1\cdots x_r D$ for every atom $a$. Since $\tau^p $ induces a
permutation on the set of atoms, this means that $a\preccurlyeq
x_1\cdots x_r D$ for every atom $a$. But since $G$ is atom-friendly,
$\Delta=\mbox{lcm}(\mathcal A)$, hence we finally obtain that
$\Delta \preccurlyeq x_1\cdots x_r D$.

Now denote $z_1\cdots z_r$ the right normal form of $x_1\cdots x_r$.
We just showed that $\Delta \preccurlyeq z_1\cdots z_r D$, but this
is equivalent to say that $z_1\cdots z_r D \succcurlyeq \Delta$.
Since $z_1\cdots z_r$ is in right normal form, this implies that
$z_r D\succcurlyeq \Delta$, which is equivalent to
$\Delta\preccurlyeq z_rD$ or, in other words, $\partial(z_r)
\preccurlyeq D$.

Now we use again that $G$ is atom-friendly, so ${\cal S}(D) =
\mathcal D$.  But since ${\cal D} = {\cal A} \backslash {\cal
S}(\iota(x)))$, one has that ${\cal S}(D) \cap {\cal
S}(\iota(x))=\emptyset$. This means that $D \wedge_L \iota(x)=
D\wedge_L \tau^{-p}(x_1) =1$, which is equivalent to $\tau^p(D)
\wedge_L x_1 =1$.

Finally, consider $y=\mathbf c_R(x) = x^{z_r^{-1}}= \Delta^p
\tau^p(z_r) z_1\cdots z_{r-1}$. We will show that $\mathbf
c_L(y)=x$. Recall that $\partial(z_r) \preccurlyeq D$, hence
$\partial(\tau^p(z_r)) \preccurlyeq \tau^p(D)$. On the other hand,
$z_1\cdots z_{r-1} \preccurlyeq z_1\cdots z_r =x_1\cdots x_r$.
Hence, if we denote by $\alpha=\iota(z_1\cdots z_{r-1})$, we have
$\alpha \preccurlyeq \iota(z_1\cdots z_r)= \iota(x_1\cdots x_r) =
x_1$.  But since $\tau^p(D)\wedge_L x_1 =1$, and we are considering
left divisors $\partial(\tau^p(z_r)) \preccurlyeq \tau^p(D)$ and
$\alpha \preccurlyeq x_1$, it follows that $\partial(\tau^p(z_r))
\wedge_L \alpha =1$. In other words, $\tau^p(z_r) \alpha$ is left
weighted as written. This is equivalent to say that $\tau^p(z_r)$ is
the first factor in the left normal form of $\tau^p(z_r) z_1\cdots
z_{r-1}$. Therefore $\mathbf c_L(y) = y^{z_r} = x$, as we wanted to
show.
\end{proof}

We have thus shown Theorem~\ref{T:surjective}, since Propositions
\ref{P:preimage1} and \ref{P:preimage2} run over all possibilities.

We end this section by recalling a result by Maffre~\cite{Maffre}
showing when each of the above two cases hold.

\begin{theorem}{\rm \cite{Maffre}}
Let $G$ be a Garside groups, and let $x=\Delta^p x_1\cdots x_r\in G$
be written in left normal form. Then
\begin{enumerate}

\item $\mathbf c_L(y)=x$ for some $y\in G$ with $\inf(y)=p-1$, if and only
if $\mathbf c_L(\tau^{-1}(x^a))=x$ \mbox{for some atom $a$.}

\item $\mathbf c_L(y)=x$ for some $y\in G$ with $\inf(y)=p$, if and only if
$\mathbf c_L(\mathbf c_R(x))=x$.

\end{enumerate}
\end{theorem}

What we showed in Theorem~\ref{T:surjective} is that at least one of
the above cases must happen.

\section{Rigid ultra summit sets}\label{S:rigid_USS}

\subsection{Left rigid and right rigid elements}

In this section we will show Theorem~\ref{T:rigid_USS}. Let $x\in
B_n$, and recall the definition of $USS_L(x)$ and $USS_R(x)$ given
in Section~\ref{S:Garside}. Since the statement of
Theorem~\ref{T:rigid_USS} refers to the conjugacy class of $x$, and
not to $x$ itself, we can assume that $x\in SSS(x)$, that is, $x$
has minimal canonical length in its conjugacy class. We will see how
one can determine if $x$ is conjugate to a rigid braid by looking at
its powers.  First we will see that if $x$ is conjugate to a rigid
element, then the infimum and supremum of its powers behave as one
should expect.

\begin{definition}{\rm \cite{Lee_transition}}  Given an element $x$ in a
Garside group $G$, we say that $x$ is {\bf periodically geodesic} if
$\inf(x^m)=m\inf(x)$ and $\sup(x^m)=m\sup(x)$ for every $m\geq 1$.
\end{definition}

\begin{lemma}\label{L:conjugate_to_rigid->pg}
If $x\in SSS(x)$ in a Garside group $G$, and $x$ is conjugate to a
(left or right) rigid element, then $x$ is periodically geodesic.
\end{lemma}

\begin{proof}
Let $y=\Delta^p y_1\cdots y_r$ be a left rigid element conjugate to
$x$. Then every power of $y$ is left rigid and $y$ is periodically
geodesic. Notice also that the left normal form of $x$ is
$x=\Delta^p x_1\cdots x_r$, where $p$ and $r$ are the same as above,
since $x\in SSS(x)$. Hence $\inf(x^m)\geq pm$ and $\sup(x^m)\leq
(p+r)m$. Now $y^m$ is rigid, thus $y^m \in USS(y^m)\subset
SSS(y^m)$, hence $\inf(y^m)=pm$ is maximal in its conjugacy class,
and $\sup(y^m)=(p+r)m$ is minimal in its conjugacy class. Since
$x^m$ is conjugate to $y^m$, this implies that $\inf(x^m)= pm = m
\inf(x)$ and $\sup(x^m)=(p+r)m= m\sup(x)$, so $x$ is periodically
geodesic.
\end{proof}

The above result is not the only one relating periodically geodesic
and rigid elements.

\begin{lemma}\label{L:pg_&_rigid_power->rigid}
Let $x$ be an element in a Garside group $G$. If $x$ is periodically
geodesic and $x^m$ is left (resp. right) rigid for some $m\geq 1$,
then $x$ is left (resp. right) rigid.
\end{lemma}

\begin{proof}
Let $\Delta^p x_1\cdots x_r$ be the left normal form of $x$. Since
$x$ is periodically geodesic, the left normal form of $x^m$ is
$\Delta^{mp}z_1\cdots z_{rm}$, where
$$
z_1\cdots z_{rm} = \tau^{(m-1)p}(x_1\cdots x_r)
\tau^{(m-2)p}(x_1\cdots x_r)\cdots \tau^p(x_1\cdots x_r)(x_1\cdots
x_r).
$$
This means that $\tau^{(m-1)p}(x_1)\preccurlyeq z_1\cdots z_{rm}$,
hence $\tau^{(m-1)p}(x_1)\preccurlyeq z_1$, since $z_1\cdots z_{rm}$
is in left normal form. But then
$\iota(x)=\tau^{-p}(x_1)\preccurlyeq  \tau^{-mp}(z_1) =\iota(x^m)$.

In the same way, since the last simple factor in the above
decomposition of $z_1\cdots z_{rm}$ is $x_r$, and the number of
factors is precisely $rm$, it follows that $x_r\succcurlyeq z_{rm}$.
In other words, $\varphi(x)\succcurlyeq \varphi(x^m)$.

% --
%Finally, recall that $x^m$ is rigid, which means that
%$\varphi(x^m)\iota(x^m)$ is left weighted as written, that is, for
%every prefix $a\preccurlyeq \iota(x^m)$, the product $\varphi(x^m)a$
%is not simple. Since $\iota(x)\preccurlyeq \iota(x^m)$, this implies
%that for every prefix $a\preccurlyeq \iota(x)$, the product
%$\varphi(x^m)a$ is not simple. But since $\varphi(x)\succcurlyeq
%\varphi(x^m)$, it follows that $\varphi(x)a$ is not simple for every
%$a\preccurlyeq \iota(x)$. That is, $\varphi(x)\iota(x)$ is left
%weighted,  hence $x$ is rigid.
Finally, recall that $x^m$ is rigid, which means that
$\varphi(x^m)\iota(x^m)$ is left weighted as written, that is,
$\partial(\varphi(x^m)) \wedge_L \iota(x^m) = 1$.  Since
$\varphi(x)\succcurlyeq \varphi(x^m)$ is equivalent to
$\partial(\varphi(x))\preccurlyeq\partial(\varphi(x^m))$, we have
$\partial(\varphi(x)) \wedge_L \iota(x) \preccurlyeq
\partial(\varphi(x^m)) \wedge_L \iota(x^m) = 1$. That is,
$\varphi(x)\iota(x)$ is left weighted, whence $x$ is rigid.
% -- VG
\end{proof}

\begin{corollary}\label{C:left_right_rigid}
% --
%Let $x$ be an element of a Garside group $G$ with $\ell(x)>1$. If
Let $x$ be an element of a Garside group $G$. If
% -- VG
$x$ has a left rigid power and $x$ is conjugate to a right rigid
element, then $x$ if left rigid. Also, if $x$ has a right rigid
power and $x$ is conjugate to a left rigid element, then $x$ is
right rigid.
\end{corollary}

\begin{proof}
This is a direct consequence of
Lemmas~\ref{L:conjugate_to_rigid->pg} and
\ref{L:pg_&_rigid_power->rigid}.
\end{proof}

After this result, in order to show that every left rigid element is
conjugate to a right rigid element, and vice versa, we must show that
every left rigid element has a conjugate which has a right rigid
power.  In braid groups, this holds for pseudo-Anosov braids, since
one has the following result.

\begin{theorem}\label{T:pA_USS->rigid power}{\rm \cite[Theorem 3.23]{BGG1}}
Let $x\in B_n$ be a pseudo-Anosov braid. If $x \in USS_L(x)$ and
$\ell(x)>1$, then $x$ has a left rigid power. In the same way, if
$x\in USS_R(x)$ and $\ell(x)>1$, then $x$ has a right rigid power.
\end{theorem}

\begin{corollary}\label{C:pA->Lrigid=Rrigid}
If $x\in B_n$ is a left (resp. right) rigid, pseudo-Anosov braid,
and $\ell(x)>1$, then $x$ is conjugate to a right (resp. left) rigid
braid.
\end{corollary}

\begin{proof}
Suppose that $x$ is left rigid, and consider $y\in USS_R(x)$. By
% --
%Theorem~\ref{T:pA_USS->rigid power} $y$ has a rigid power, hence by
%Corollary~\ref{C:left_right_rigid} $y$ itself must be rigid.  If $x$
Theorem~\ref{T:pA_USS->rigid power}, the braid $y$ has a right rigid
power, hence $y$ itself must be right rigid by
Corollary~\ref{C:left_right_rigid}.  If $x$
% -- VG
is right rigid, the proof follows the same reasoning.
\end{proof}

But there are two more kind of braids, namely periodic and reducible
ones. Does the above result hold for these ones? The answer is
positive, as we shall see.  We recall that a braid $x\in B_n$ is
called {\bf periodic} if $x^m=\Delta^t$ for some nonzero integers
$m$ and $t$. The above result holds trivially for periodic braids,
due to the following lemma.

\begin{lemma}\label{L:periodic->nonrigid}
A left or right rigid braid can never be periodic.
\end{lemma}

\begin{proof}
By definition, if $x\in B_n$ is rigid then $\ell(x)>0$. Also, by
Lemma~\ref{L:conjugate_to_rigid->pg}, $x$ is periodically geodesic.
Hence $\ell(x^m)=|m|\ell(x)>0$ for every nonzero integer $m$.
Therefore no power of $x$ can be a power of $\Delta$, since
$\ell(\Delta^t)=0$ for every $t$.
\end{proof}

It just remains to show the case of reducible braids. A braid $x\in
B_n$ is said to be {\bf reducible} if, regarding $x$ as a
homeomorphism of the $n$-times punctured disc, it preserves a family
of disjoint, closed, essential curves, up to isotopy~\cite{BLM}.
This can be expressed in other terms: A braid $x\in B_n$ is said to
admit a {\bf coherent tape structure}~\cite{BGN} if it can be
obtained from a braid $\widehat x\in B_m$, with $m<n$, by replacing,
for each $i=1,\ldots,m$, the $i$-th strand of $\widehat x$ by a
braid $x_{[i]}\in B_{k_i}$, with $k_i\geq 1$. One can think that the
$i$-th strand of $\widehat x$ becomes a tube, and that $x_{[i]}$
lies inside that tube. One further requirement is that the $m$-tuple
$(k_1,\ldots, k_m)$ is coherent with the permutation induced by
$\widehat x$, that is, if the $i$-th strand of $\widehat x$ ends at
position $j$, then $k_i=k_j$. The braid $\widehat x$ is called the
{\bf tubular}, or {\bf external} braid of this decomposition of $x$,
while each $x_{[i]}$ is called the $i$-th {\bf internal} braid. A
braid is then {\bf periodic} if one of its conjugates admits a
coherent tape structure.

We can now extend the result of Corollary~\ref{C:pA->Lrigid=Rrigid}
to the whole $B_n$, so we can show the following result, which is equivalent to Theorem~\ref{T:rigid_USS}.

\begin{theorem}\label{T:left_right_rigid}
If $x\in B_n$ is a left (resp. right) rigid braid, and $\ell(x)>1$,
then $x$ is conjugate to a right (resp. left) rigid braid.
\end{theorem}

\begin{proof}
Suppose that $x$ is left rigid. We will show the result by induction
on $n$. If $n=1$, $x$ is trivial and there is nothing to show. If
$n=2$, $x$ is either trivial or periodic and by
Lemma~\ref{L:periodic->nonrigid}, it cannot be rigid. We then
suppose that $n>2$ and that the result holds for braids with less
than $n$ strands.

If $x$ is pseudo-Anosov, the result is given by
Corollary~\ref{C:pA->Lrigid=Rrigid}. On the other hand, $x$ cannot
be periodic by Lemma~\ref{L:periodic->nonrigid}. Hence we can assume
that $x$ is reducible.

In~\cite{BGN} it was shown that if a braid $\alpha$ admits a
coherent tape structure, so do $\mathbf c_L(\alpha)$ and $\mathbf
d_L(\alpha)$. By symmetry, the same property holds for $\mathbf
c_R(\alpha)$ and $\mathbf d_R(\alpha)$. This implies that for every
reducible braid, there is some element in its (left or right) ultra
summit set that admits a coherent tape structure. Since we are
assuming that $x$ is left rigid and $\ell(x)>1$,  $USS_L(x)$ is
the set of left rigid conjugates of $x$, hence there is a conjugate
of $x$ which is left rigid, and admits a coherent tape structure. We
can then assume that $x$ itself admits a coherent tape structure.

Let $y\in USS_R(x)$, obtained from $y$ by a finite number of
applications of $\mathbf c_R$ and $\mathbf d_R$. After~\cite{BGN},
$y$ admits a coherent tape structure. By
Corollary~\ref{C:left_right_rigid}, we just need to show that $y$
has a right rigid power.

We will denote $\widehat y\in B_m$ and $y_{[1]},\ldots, y_{[m]}$,
respectively, the external and internal braids associated to $y$, where
$y_{[i]}\in B_{k_i}$ for $i=1,\ldots, m$, and $k_1+\cdots +k_m=n$.
Notice that if $y$ admits a coherent tape structure, so does every
power of $y$. In order to simplify the notation, we will replace $y$
by a power $y^m$ such that the permutation induced by
$\widehat{y^{m}}$ is trivial ($\widehat{y^{m}}$ is a pure braid).
Notice that $x^m$ is left-rigid, $y^m$ admits a coherent tape
structure, and if we show that $y^m$ has a right rigid power, this
will also be true for $y$. Hence can assume that $\widehat y$ is a
pure braid.

Let $p=\inf(x)$ and $p+r=\sup(x)>1$. Notice that, since $x$ is left
rigid, $\varphi(x)\iota(x)$ is left weighted. One can see the tape
structure of $x$ in this pair of simple elements, in the following
way. One has $\inf(x)=\min\{\inf(\widehat
x),\inf(x_{[1]}),\ldots,\inf(x_{[m]})\}$ and
$\sup(x)=\max\{\sup(\widehat
x),\sup(x_{[1]}),\ldots,\sup(x_{[m]})\}$. The part of $\widehat x$
(resp. $x_{[i]}$) that one can see in $\varphi(x)$ will be
$\varphi(\widehat x)$ (resp. $\varphi(x_{[i]})$) if $\sup(\widehat
x)=p+r$ (resp. $\sup(x_{[i]})=p+r$), and will be trivial otherwise.
Analogously, the part of $\widehat x$ (resp. $x_{[i]}$) that one can
see in $\iota(x)$ will be $\iota(\widehat x)$ (resp.
$\iota(x_{[i]})$) if $\inf(\widehat x)=p$ (resp. $\inf(x_{[i]})=p$),
and will be equal to $\Delta\in B_m$ (resp. $\Delta\in B_{k_i}$)
otherwise. If we had a trivial component in $\varphi(x)$, then
$\varphi(x)\iota(x)$ could not be left weighted, unless the
corresponding component of $\iota(x)$ would be trivial. In the same
way, If we had a  $\Delta$ component in $\iota(x)$, then
$\varphi(x)\iota(x)$ could not be left weighted, unless the
corresponding component of $\iota(x)$ would be also $\Delta$.
Therefore, each external or internal component of $x$ must be as
follows: either it is trivial, or it is $\Delta^{p+r}$ (with the
corresponding number of strands), or it is left rigid with infimum
$p$ and supremum $p+r$. This has an important consequence: applying
(left or right) cyclings and decyclings to $x$ induces (left or
right) cyclings and decyclings to $\widehat x$, $x_{[1]},\ldots,
x_{[m]}$.   Therefore, $y\in USS_R(x)$ implies that $\widehat y\in
USS_R(\widehat x)$ and $y_{[i]}\in USS_R(x_{[i]})$ for
$i=1,\ldots,m$.

Finally, each of the components of $y$ (having less than $n$
strands) which is neither trivial nor $\Delta^{p+r}$ is conjugate to
a left rigid braid with canonical length greater than 1. The
induction hypothesis tells us that each of these components is then
right rigid, and it has infimum $p$ and supremum $p+r$. Therefore,
$y$ itself must be right rigid, as we wanted to show.
\end{proof}

\subsection{Left and right ultra summit graphs are isomorphic}

We will now show that given a left rigid braid $x\in USS_L(x)$ with
$\ell(x)>1$, then the directed graphs $USG_L(x)$
and $USG_R(x)$ are isomorphic, with the arrows reversed. That is, we
will show Theorem~\ref{T:rigid_USG}. We need to define an isomorphism
of directed graphs (in other words, an invertible functor from the
category $USG_L(x)$ to the category $USG_R(x)^{op}$). The isomorphism
is very easy to define at the level of vertices (objects), that is,
the elements of the ultra summit sets.

% --  changes reverted for now
\begin{definition}
Let $x\in B_n$ be a left rigid braid, with $\ell(x)=r>1$. We define
$\Phi(x)=\mathbf c_R^{2rt}(x)$, where $t$ is any non-negative integer
such that $\mathbf c_R^{2rt}(x)$ is right rigid.
\end{definition}
%In fact, all arguments in this section are valid for the case of a
%general Garside group, under the additional assumption that $x$ has a
%right rigid conjugate.  (We know by Theorem~\ref{T:rigid_USS} that
%this additional assumption is not needed for braids.)  Recall that
%there is a smallest positive integer $e$ such that $\Delta^e$ is
%central; in the case of $B_n$, we have $e=2$.
%
%\begin{definition}
%Let $G$ be a Garside group and let $x\in G$ be left rigid braid, which
%has a right rigid conjugate and satisfies $\ell(x)=r>1$. We define
%$\Phi(x)=\mathbf c_R^{ert}(x)$, where $t$ is any non-negative integer
%such that $\mathbf c_R^{ert}(x)$ is right rigid.
%\end{definition}
% -- VG

Notice that $\Phi$ is well defined: Since $x$ is left rigid, $x\in
SSS(x)$, so one can go from $x$ to $USS_R(x)$ by iterated right
cycling. Since $\ell(x)>1$, Theorem~\ref{T:left_right_rigid} tells
us that $x$ is conjugate to a right rigid element, hence $USS_R(x)$
consists of right rigid elements, and one obtains a right rigid
element by applying iterated right cycling to $x$. Also, for every
right rigid element $z$ with $\ell(z)=r$, one has $\mathbf
c_R^{2r}(z)=z$. Hence, if $t$ is an integer such that $\mathbf
c_R^{2rt}(x)$ is right rigid, then $\mathbf c^{2rt}_R(x)= \mathbf
c^{2r}(\mathbf c^{2rt}_R(x))=\mathbf c^{2r(t+1)}_R(x)$. This implies
that if $\mathbf c^{2rt}_R(x)$ and $\mathbf c^{2rt'}_R(x)$ are both
right rigid, they are equal. Hence $\Phi$ is well defined.
% --
%Notice that in an arbitrary Garside group, $2rt$ should be replaced by
%$ert$, where $e$ is the smallest integer such that $\Delta^e$ is
%central.
% -- VG

We will show below that $\Phi$ is a bijective map from $USS_L(x)$ to
$USS_R(x)$. But we also want to show that $USG_L(x)$ is isomorphic to $USG_R(x)^{op}$.
We already know a map $\Phi$ that sends vertices (objects) of $USG_L(x)$ to vertices (objects) of $USG_R(x)^{op}$.  Let us see how $\Phi$ is defined on the arrows
(morphisms) of $USG_L(x)$.  In order to do this, we recall the
% --
%definition of Gebhardt's transport. This map is defined
definition of the transport map. This map is defined
% -- VG
in~\cite{Gebhardt} using left normal forms, but it can be equally
defined, by symmetry, using right normal forms.

\begin{definition}{\rm \cite{Gebhardt}}
Given $x\in SSS(x)$ in a Garside group, and given a positive element
$u$ such that $u^{-1}xu=y \in SSS(x)$, one defines the {\bf left
transport} of $u$ as:
$$
 u^{(1)}_L = \iota_L(x)^{-1}\cdot u \cdot
\iota_L(y).
$$
The {\bf iterated left transports} of $u$ are defined recursively,
for every $i\geq 1$, by
$$
u^{(i)}_L=\left(u^{(i-1)}_L\right)_L^{(1)}.
$$
\end{definition}

Notice that, since $u^{-1}xu=y$, one has
$\left(u^{(i)}_L\right)^{-1}\mathbf c_L^i(x)\; u^{(i)}_L=\mathbf
c_L^i(y)$. In other words, since $u$ conjugates  (on the
right) $x$ to $y$ , the $i$-th left transport of $u$ conjugates (on the right) the $i$-th left cycling of $x$ to the $i$-th left cycling of $y$.

\begin{definition}{\rm \cite{Gebhardt}}
Given $x\in SSS(x)$ in a Garside group, and given a positive element
$v$ such that $vxv^{-1}=z \in SSS(x)$, one defines the {\bf right
transport} of $v$ as:
$$
 v^{(1)}_R = \iota_R(z)\cdot v \cdot
\iota_R(x)^{-1}.
$$
The {\bf iterated right transports} of $v$ are defined recursively,
for every $i\geq 1$, by
$$
v^{(i)}_R=\left(v^{(i-1)}_R\right)_R^{(1)}.
$$
\end{definition}

In this case, since $vxv^{-1}=z$, one has $v^{(i)}_R \: \mathbf
c_R^i(x)\left(v^{(i)}_R\right)^{-1}=\mathbf c_R^i(z)$. In other
words, since $v$ conjugates (on the left) $x$ to $z$, the $i$-th
right transport of $v$ conjugates (on the left) the $i$-th right cycling of $x$ to
the $i$-th right cycling of $z$.

\begin{theorem}{\rm \cite{Gebhardt}}
With the above conditions, one has the following properties, for
every $i\geq 1$:

\begin{tabular}{ll}
  {\it 1.}  If $u_1\preccurlyeq u_2$ then $(u_1)^{(i)}_L\preccurlyeq (u_2)^{(i)}_L$. \hspace{1cm} \vspace{.2cm}
 & If $v_1\succcurlyeq v_2$ then $(v_1)^{(i)}_R\succcurlyeq (v_2)^{(i)}_R$.
 \\ \vspace{.2cm} {\it 2.}  $\left(u_1\wedge_L u_2\right)_L^{(i)} = (u_1)^{(i)}_L \wedge_L
 (u_2)_L^{(i)}$.  &  $\left(v_1\wedge_R v_2\right)_R^{(i)} = (v_1)_R^{(i)} \wedge_R  (v_2)_R^{(i)}$.
 \\ \vspace{.2cm} {\it 3.} $\Delta^{(i)}_L = \Delta,   \hspace{1cm}   1^{(i)}_L = 1.$  &  $\Delta^{(i)}_R = \Delta,   \hspace{1cm}  1^{(i)}_R=1$.
 \\ {\it 4.}   If $u$ is simple, $u^{(i)}_L$ is simple.  &  If $v$ is simple, $v^{(i)}_R$ is simple.
\end{tabular}
\end{theorem}

Let us then define $\Phi$ on the arrows of $USS_L(x)$.

\begin{definition}
Let $x,y\in USS_L(x)\subset B_n$ be left rigid braids with
$\ell(x)>1$, and let $t$ be a nonnegative integer such that
$\Phi(x)=\mathbf c_R^{2rt}(x)$ and $\Phi(y)=\mathbf c_R^{2rt}(y)$.
Given $u\in B_n$ such that $u^{-1}xu= y$, so that $uyu^{-1}=x$, we
define $\Phi(u)= u^{(2rt)}_R$.
\end{definition}

\begin{proposition}
$\Phi$ is a well defined map of directed graphs (a well defined functor) from
$USG_L(x)$ to $USG_R(x)^{op}$.
\end{proposition}

\begin{proof}
We already know that $\Phi(y)\in USS_R(x)$ for every $y\in
USS_L(x)$, hence $\Phi$ sends vertices of $USG_L(x)$ to vertices of
$USG_R(x)^{op}$. Now consider an arrow $s$ going from $x$ to $y$ in
$USG_L(x)$. Since $s^{-1}xs=y$, one has $sys^{-1}=x$. Hence, if we
denote $s_0=s^{(2rt)}_R$ for an integer $t$  such that $\Phi(x)=\mathbf c_R^{2rt}(x)$
and $\Phi(y)=\mathbf c_R^{2rt}(y)$, we have $s_0 \:\mathbf
c^{2rt}_R(y)\: s_0^{-1}= \mathbf c_R^{2rt}(x)$, that is, $s_0 \:\Phi(y)\:
s_0^{-1} = \Phi(x)$, where $\Phi(y)$ and $\Phi(x)$ are right rigid.

Notice that, since $\Phi(y)$ is right rigid and has canonical length
$r$, then $\mathbf c_R^{2r}(\Phi(y))=\Phi(y)$, since the product of
the $2r$ conjugating elements for right cycling is precisely
$\Phi(y)^2\Delta^{-2}$. In the same way, the product of the $2r$ conjugating
elements that perform iterated right cycling of $\Phi(x)$ is
precisely $\Phi(x)^2\Delta^{-2}$. Hence, the $2r$-th iterated right transport of
$s_0$ is $s_0^{(2r)}= \Phi(x)^2\Delta^{-2} s_0 \Delta^2\Phi(y)^{-2} = \Phi(x)^2 s_0 \Phi(y)^{-2} = \Phi(x) s_0 \Phi(y)^{-1} = s_0$. This means
that $s^{(2rt')}=s^{(2rt)}$ for every $t'\geq t$. Hence $\Phi(s)$ is
a well defined simple element which is, by the above argument, an
arrow in $USG_R(x)$ going from $\Phi(y)$ to $\Phi(x)$, hence an arrow in $USG_R(x)^{op}$ going from $\Phi(x)$ to
$\Phi(y)$.
\end{proof}

It remains to show that $\Phi$ is invertible.  In order to do this,
we start by recalling a result from~\cite{BGG1} that relies cyclings
and powers. Given $x$ in a Garside group $G$, denote
$C_i=\iota(\mathbf c_L^{i-1}(x))$ for every $i\geq 1$. That is, $C_i$
is the conjugating element from $\mathbf c_L^{i-1}(x)$ to $\mathbf
c_L^i(x)$, and $x^{C_1\cdots C_i} = \mathbf c_L^i(x)$. Then one has:

\begin{lemma}{\rm \cite[Lemma 2.4]{BGG1}}\label{L:power_LNF} Let $G$ be a Garside group and let
$ x\in SSS(x)\subset G$, with $\inf(x)=p$ and $\ell(x)>1$. Then, for
every $m\geq 1$,
$$
     x^m \Delta^{-mp}= C_1\cdots C_m  \mathbf R_m,
$$
where \begin{enumerate}
 \item  $\sup(C_1\cdots C_m)=m$ and $\varphi_L(C_1\cdots C_m)\succcurlyeq \varphi_L(\mathbf c_L^m(x))$.

 \item $\inf(\mathbf R_m)=0$ and $\iota_L(\mathbf R_m)\preccurlyeq
 C_{m+1}=\iota_L(\mathbf c_L^m(x))$.
\end{enumerate}
\end{lemma}

This result can be improved if $x$ is conjugate to a rigid element.

\begin{lemma}\label{L:conjugate_to_rigid}
Let $G$ be a Garside group and let $ x\in SSS(x)\subset G$, with
$\inf(x)=p$ and $\ell(x)>1$. Suppose that $x$ is conjugate to a left
rigid element, and let $m$ be such that $y=\mathbf c_L^m(x)$ is
rigid. Then
$$
    C_1\cdots C_m = (x^m \Delta^{-mp})\wedge_L \Delta^m,
$$
where $\inf(C_1\cdots C_m)=0$ and $\sup(C_1\cdots C_m)=m$.
\end{lemma}

\begin{proof}
By the above lemma, $C_1\cdots C_m \preccurlyeq  x^m \Delta^{-mp}$.
But since $m$ is conjugate to a rigid element,
Lemma~\ref{L:conjugate_to_rigid->pg} implies that $\inf(x^m)=mp$, so
$\inf(x^m\Delta^{-mp})=0$. This means that $\inf(C_1\cdots C_m)=0$.

Recall also that $x^m \Delta^{-pm}= C_1\cdots C_m \mathbf R_m $,
where $\varphi_L(C_1\cdots C_m)\succcurlyeq \varphi_L(\mathbf
c^m(x))=\varphi_L(y)$ and $\iota_L(\mathbf R_m)\preccurlyeq
\iota_L(\mathbf c^m(x)) = \iota_L(y)$. Since $y$ is left rigid,
the decomposition $\varphi_L(y) \iota_L(y)$ is left weighted. Hence,
if $z_1\cdots z_m$ is the left normal form of $C_1\cdots C_m$, this
means that $z_1\cdots z_m \iota_L(\mathbf R_m)$ is in left normal
form as written. In other words, the first $m$ factors of the left
normal form of $x^m\Delta^{-mp}$ are precisely $z_1\cdots z_m
=C_1\cdots C_m$. That is, $C_1\cdots C_m = (x^m \Delta^{-mp})
\wedge_L \Delta^m$, as we wanted to show.
\end{proof}

This allows us to determine very precisely the left normal form of
$x^m$, for $m$ big enough, when $x$ is conjugate to a left rigid
element.  In order to avoid confusing notation produced by the
powers of $\Delta$ in the normal forms, we will introduce the
following notion:

\begin{definition}
% --
%Given an element $z\in B_n$, whose left normal form is
Let $G$ be a Garside group.
Given an element $z\in G$, whose left normal form is
% -- VG
$\Delta^pz_1\cdots z_r$ and whose right normal form is
$z_1'\cdots z_r'\Delta^p$, we define the {\bf left interior} of $z$ as
$$
z^\circ_L \quad  =  \quad z \Delta^{-p}\quad  = \quad
\tau^{-p}(z_1)\cdots \tau^{-p}(z_r) \quad = \quad z_1'\cdots z_r',
$$
and the {\bf right interior} of $z$ as
$$
z^\circ_R  \quad =  \quad  \Delta^{-p} z \quad  = \quad z_1\cdots
z_r \quad =\quad \tau^{p}(z_1')\cdots \tau^{p}(z_r').
$$
\end{definition}

Notice that the above factorizations are, respectively, the left and
right normal forms of $z_L^{\circ}$ and of $z_R^{\circ}$. Notice
also that if $y= \Delta^p y_1\cdots y_r$ is left rigid, then
$$
(y^m)^\circ_L = y^m \Delta^{-pm}= \left(\tau^{-p}(y_1)\cdots
\tau^{-p}(y_r)\right)\left(\tau^{-2p}(y_1)\cdots
\tau^{-2p}(y_r)\right)\cdots \left(\tau^{-mp}(y_1)\cdots
\tau^{-mp}(y_r)\right),
$$
and it is in left normal form as written. Moreover, in this case
$(y^m)^\circ_L$ is precisely the conjugating element that takes $y$
to $\mathbf c^{rm}_L(y)$.

\begin{lemma}\label{L:conjugate_to_rigid_LNF}
Let $G$ be a Garside group and let $ x\in SSS(x)\subset G$, with
$\inf(x)=p$ and $\ell(x)=r>1$. Suppose that $x$ is conjugate to a
left rigid element. Let $N$ be such that $y=\mathbf c_L^N(x)$ is
left rigid. Then:
\begin{enumerate}

\item There exists an integer $M$ such that $\left(y^M\right)^\circ_R \succcurlyeq
C_1\cdots C_N$.

\item Let $M$ be an integer satisfying the above condition. If $z_1\cdots z_N$ is the left normal form of $C_1\cdots C_N$, and
$z_1'\cdots z_s'$ is the left normal form of
$\left(y^M\right)^\circ_R (C_1\cdots C_N)^{-1}$, then for every
$m\geq M$, the left normal form of $(x^m)^{\circ}_L$ is
$$
  (x^m)^{\circ}_L =  \left(z_1\cdots z_N \right)\cdot
  \left(y^{m-M}\right)^\circ_L \cdot \left( \tau^{-pm}(z_1')\cdots \tau^{-pm}(z_s')\right),
$$
where the central factor is assumed to be written in left normal form. Moreover, $N+s=Mr$.
\end{enumerate}
\end{lemma}

\begin{proof}
Recall that $x^{C_1\cdots C_N}= \mathbf c_L^N(x)=y$, so
$\left(x^N\right)^{C_1\cdots C_N}=y^N$. Recall also by
Lemma~\ref{L:conjugate_to_rigid} that $C_1\cdots C_N \preccurlyeq
\left(x^N\right)^\circ_L = x^N \Delta^{-pN}$. This means that
$\alpha=(C_1\cdots C_N)^{-1} \left(x^N\right)^\circ_L $ is a
positive braid. Hence $ y^N= (C_1\cdots C_N)^{-1} x^N (C_1\cdots
C_N) = \alpha \Delta^{pN} C_1\cdots C_N= \Delta^{pN}
\tau^{pN}(\alpha) C_1\cdots C_n, $ so $\left(y^N\right)^\circ_R =
\Delta^{-pN}y^N = \tau^{pN}(\alpha) C_1\cdots C_n \succcurlyeq
C_1\cdots C_N$. Hence the first property is satisfied for $M=N$.

Now let $M$, $m$, $z_1\cdots z_N$ and $z_1'\cdots z_s'$ be defined
as in Condition 2. Notice that since $m\geq M$, one has
$\left(y^m\right)^\circ_R = \Delta^{-pm}y^m \succcurlyeq
\Delta^{-pM}y^M \succcurlyeq C_1\cdots C_N$. That is, there exists a
positive braid $\beta$ such that $y^m = \Delta^{mp} \beta C_1\cdots
C_N$.  Since $y$ is a left rigid element,  by Lemma~\ref{L:power_LNF}, $\varphi_L(C_1\cdots C_N)\succcurlyeq \varphi_L(\mathbf c_L^N(x))=\varphi_L(y)$.
Also, $\iota(\tau^{-mp}(\beta))\preccurlyeq \iota(y^m)=\iota(y)$. This
implies, as $\varphi(y)\iota(y)$ is left weighted, that $z_N
\iota(\tau^{-mp}(\beta))$ is also left weighted.

If we now conjugate $y^m$ by $(C_1\cdots C_N)^{-1}$, we obtain $x^m=
C_1\cdots C_N \Delta^{mp} \beta= C_1\cdots C_N \tau^{-mp}(\beta)
\Delta^{mp}$, hence $\left(x^m\right)^\circ = C_1\cdots C_N
\tau^{-mp}(\beta)= z_1\cdots z_N \tau^{-mp}(\beta)$. Since $z_N
\iota(\tau^{-mp}(\beta))$ is left weighted, it follows that the
first $N$ factors in the left normal form of
$\left(x^m\right)^\circ_L$ are precisely $z_1\cdots z_N$.

Now recall that $z_1'\cdots z_s'$ is the left normal form of
$\Delta^{-pM} y^M (C_1\cdots C_N)^{-1}$. Hence
$$
y^m= y^{m-M} y^M = y^{m-M}\Delta^{pM} z_1'\cdots z_s' C_1\cdots C_N =
\left(y^{m-M}\right)^\circ_L \Delta^{p(m-M)} \Delta^{pM} z_1'\cdots z_s' C_1\cdots C_N
$$
$$
  = \left(y^{m-M}\right)^\circ_L \Delta^{pm} z_1'\cdots z_s' C_1\cdots C_N
  = \left(y^{m-M}\right)^\circ_L \left( \tau^{-pm}(z_1')\cdots \tau^{-pm}(z_s')\right) \Delta^{pm} C_1\cdots C_N.
$$
Conjugating by $(C_1\cdots C_N)^{-1}$, one obtains
$$
 x^m = \left(C_1\cdots C_N \right) \left(y^{m-M}\right)^\circ_L \left( \tau^{-pm}(z_1')\cdots \tau^{-pm}(z_s')\right) \Delta^{pm},
$$
hence
$$
   (x^m)^{\circ}_L =  \left(z_1\cdots z_N \right)\cdot
  \left(y^{m-M}\right)^\circ_L \cdot \left( \tau^{-pm}(z_1')\cdots \tau^{-pm}(z_s')\right).
$$
This is written in left normal form since
$\varphi\left(\left(y^{m-M}\right)^\circ_L\right) \tau^{-pm}(z_1')$ is
left weighted, as can be seen by noticing that
$\varphi\left(\left(y^{m-M}\right)^\circ_L\right) =
\varphi\left(\tau^{-p(m-M)}(y)\right)$, and also that
$z_1'=\iota\left(\Delta^{-pM} y^M (C_1\cdots
C_N)^{-1}\right)\preccurlyeq \iota\left(\tau^{pM}(y)\right)$, so
$\tau^{-pm}(z_1')\preccurlyeq \iota\left(\tau^{-p(m-M)}(y)\right)$.

Finally, since $y$ is left rigid, $x$ is periodically geodesic. Hence $\ell(x^m)=\ell\left((x^m)^\circ_L\right)=mr$. But we just computed the left normal form of $(x^m)^\circ_L$, which has $N+(m-M)r+s$ factors. Therefore $N+(m-M)r+s=mr$, so $N+s=Mr$, as we wanted to show.
\end{proof}

By symmetry, one has the analogous result for conjugates of right rigid braids, but we will perform a slight modification:

\begin{lemma}\label{L:conjugate_to_rigid_RNF}
Let $G$ be a Garside group and let $ x\in SSS(x)\subset G$, with
$\inf(x)=p$ and $\ell(x)=r>1$. Suppose that $x$ is conjugate to a
right rigid element. Let $N$ be such that $y=\mathbf c_R^N(x)$ is
right rigid, and let $C_1',\cdots,C_N'$ the conjugating elements for
the $N$ right cyclings, that is, $(C_N'\cdots C_1') \: x \:
(C_N'\cdots C_1')^{-1} = y$. Then:
\begin{enumerate}

\item There exists an integer $M$ such that $C_N'\cdots C_1' \preccurlyeq
\left(y^M\right)^\circ_L $.

\item Let $M$ be an {\bf even} integer satisfying the above condition. If
$z_N'\cdots z_1'$ is the right normal form of $C_N'\cdots C_1'$, and
$z_s\cdots z_1$ is the right normal form of $(C_N'\cdots
C_1')^{-1}\left(y^M\right)^\circ_L $, then for every $m\geq M$, the
right normal form of $(x^m)^{\circ}_L$ is
$$
  (x^m)^{\circ}_L =  \left(z_s\cdots z_1 \right)\cdot
  \left(y^{m-M}\right)^\circ_L \cdot \left( \tau^{-pm}(z_N')\cdots \tau^{-pm}(z_1')\right),
$$
where the central factor is assumed to be written in right normal form. Moreover, \mbox{$N+s=Mr$}.
\end{enumerate}
\end{lemma}

\begin{proof}
% --
%It one applies Lemma~\ref{L:conjugate_to_rigid_LNF} to
%$\overleftarrow x$, one obtains that the right normal form of
If one follows the argument of Lemma~\ref{L:conjugate_to_rigid_LNF}
for right normal forms, one obtains that the right normal form of
% -- VG
$(x^m)^{\circ}_R$ is
$$ (x^m)^{\circ}_R =  \left(\tau^{pm}(z_s)\cdots \tau^{pm}(z_1) \right)\cdot
  \left(y^{m-M}\right)^\circ_R \cdot \left( z_N'\cdots z_1'\right),
$$
and now one just needs to notice that $(x^m)^{\circ}_L
=\tau^{-mp}\left((x^m)^{\circ}_R\right)$ and that, since $M$ is
even, $\tau^{-pm}\left(\left(y^{m-M}\right)^\circ_R\right) =
\tau^{-p(m-M)}\left(\left(y^{m-M}\right)^\circ_R\right)=
\left(y^{m-M}\right)^\circ_L$.
\end{proof}

We can now show that $\Phi$ is a bijective map on the vertices.

\begin{proposition}\label{P:Phi_bijective_on_vertices}
Let $x\in B_n$ be a left rigid braid with $\ell(x)>1$. The map
\mbox{$\Phi:\: USS_L(x) \rightarrow USS_R(x)$} defined above is
bijective.
\end{proposition}

\begin{proof}
Let us define the map $\Psi:\: USS_R(x)\rightarrow USS_L(x)$, which
is defined just as $\Phi$, by symmetry. That is,
$\Psi(z)=\overleftarrow{\Phi(\overleftarrow{z})}$. We will show that
$\Psi$ is the inverse of $\Phi$.

Let $\Delta^p x_1\cdots x_r$ be the left normal form of $x$.  Recall
that $\Phi(x)=\mathbf c_R^{2rt}(x)$ for some $t$, and then
$\Phi(x)=\mathbf c_R^{2rt'}(x)$ for every $t'\geq t$. We also have
$\Psi(\Phi(x))=\mathbf c_L^{2rs}(\Phi(x))$ for some $s$, and then
$\Psi(\Phi(x))=\mathbf c_L^{2rs'}(\Phi(x))$ for every $s'\geq s$.
Hence, if we denote $N=2r\max(t,s)$, we have $\Phi(x)=\mathbf
c_R^N(x)$ and $\Psi(\Phi(x))=\mathbf c_L^N(\Phi(x)) = \mathbf
c_L^N(\mathbf c_R^N(x))$.   We must then show that $\mathbf
c_L^N(\mathbf c_R^N(x)) = x$.

In order to do it, we will study some decompositions of $x^m$, for
$m$ big enough. For simplicity, we will consider $m$ to be even.
First, since $x$ is left rigid, the left normal form of
$\left(x^m\right)^\circ_L$ for every even $m$ is precisely:
$$
\left(x^m\right)^\circ_L=\left(\tau^{-p}(x_1) \cdots \tau^{-p}(x_r)
\right)  \left(\tau^{-2p}(x_1) \cdots \tau^{-2p}(x_r) \right) \cdots
\left(\tau^{-mp}(x_1) \cdots \tau^{-mp}(x_r) \right)
$$
$$
=  \left(\tau^{-p}(x_1) \cdots \tau^{-p}(x_r) \; x_1\cdots
x_r\right)^{m/2}
$$
$$
 = \left(\left(x^2\right)^{\circ}_L\right)^{m/2}.
$$
Notice that if $p$ is even, the above expression is just $(x_1\cdots
x_r)^m$, but if $p$ is odd this does not happen in general.

Now $x$ is conjugate to a right rigid braid, $y=\Phi(x)$. We can
then apply Lemma~\ref{L:conjugate_to_rigid_RNF} to $x$. We fix $M$
as in Lemma~\ref{L:conjugate_to_rigid_RNF}, where we can assume that
$M$ is even (otherwise, take $M+1$). We take $m$ big enough, so that
$m> 2M$ and $m$ is even. We then obtain that the right normal
form of $\left(x^m\right)^\circ_L $ is:
$$
  \left(x^m\right)^{\circ}_L =  \left(z_s\cdots z_1 \right)\cdot
  \left(y^{m-M}\right)^\circ_L \cdot \left( \tau^{-pm}(z_N')\cdots \tau^{-pm}(z_1')\right)
$$
$$
   = \left(z_s\cdots z_1 \right)\cdot
  \left(y^{m-M}\right)^\circ_L \cdot \left( z_N'\cdots z_1'\right)
$$
Notice that, by definition, $(z_N'\cdots z_1')(z_s\cdots
z_1)=(y^M)^\circ_L = y^M\Delta^{-pM}$. Also, by definition
$z_N',\ldots,z_1'$ are the conjugate elements of the iterated right
cyclings from $x$ to $y$, that is, $(z_N'\cdots z_1')\: x
\:(z_N'\cdots z_1')^{-1}= y$. Hence $\left(x^M\right)^\circ_L  = x^M
\Delta^{-pM} =  (z_N'\cdots z_1')^{-1} y^M \Delta^{-pM} (z_N'\cdots
z_1') =  (z_s\cdots z_1) (z_N'\cdots z_1')$.  Notice that we used
that $M$ is even, so $\Delta^{pM}$ is central.

We then obtain the following decomposition:
$$
(x^m)^{\circ}_L =  \left(z_s\cdots z_1 \right) \left( z_N'\cdots
z_1'\right)   \cdot   \left(x^{m-2M}\right)^\circ_L \cdot
\left(z_s\cdots z_1 \right)\left( z_N'\cdots z_1'\right).
$$
Hence
$$
   \left(y^{m-M}\right)^\circ_L = \left( z_N'\cdots
z_1'\right)   \cdot   \left(x^{m-2M}\right)^\circ_L \cdot
\left(z_s\cdots z_1 \right).
$$
Let us write the above factors in left normal form. Let $w_1\cdots
w_N$ the left normal form of $z_N'\cdots z_1'$, and let $w_1'\cdots
w_s'$ be the left normal form of $z_s\cdots z_1$. Then
$$
   \left(y^{m-M}\right)^\circ_L = \left( w_1\cdots
w_N\right)   \cdot   \left(x^{m-2M}\right)^\circ_L \cdot
\left(w_1'\cdots w_s' \right).
$$
We will now show that this decomposition is precisely the left
normal form of $\left(y^{m-M}\right)^\circ_L$.  Indeed, since
$\left(x^M\right)^\circ_L  =  (z_s\cdots z_1) (z_N'\cdots z_1') =
\left(w_1'\cdots w_s' \right)\left( w_1\cdots w_N\right)$ and
$s+N=Mr$ by Lemma~\ref{L:conjugate_to_rigid_RNF}, it follows that
the final factor of the left normal form of
$\left(x^M\right)^\circ_L$ is a suffix of $w_N$. That is,
$w_N\succcurlyeq x_r$. Since $x$ is left rigid, this implies that
$w_N \cdot \tau^{-p}(x_1)$ is left weighted, where the second factor
in this expression is the initial factor in the left normal form of
$\left(x^{m-2M}\right)^\circ_L$.  But also $w_1'$ must be a prefix of
the initial factor of $\left(x^M\right)^\circ_L$, that is,
$w_1'\preccurlyeq \tau^{-p}(x_1)$. This implies that $x_r\cdot w_1'$
is left weighted, where $x_r$ is the final factor in the left normal
form of $\left(x^{m-2M}\right)^\circ_L$.   Hence, the above
expression is the left normal form of
$\left(y^{m-M}\right)^\circ_L$, for $m$ big enough.

But recall from Lemma~\ref{L:conjugate_to_rigid} that the product of
the first $m-M$ factors in the left normal form of
$\left(y^{m-M}\right)^\circ_L$ is precisely the product of the $m-M$
conjugating elements for iterated left cycling of $y$. If we take
$m$ big enough so that $m-M\geq N$ and $m-M$ (as well as $N$) is a multiple of $2r$, the first $m-M$ factors in the left normal form of
$\left(y^{m-M}\right)^\circ_L$ are precisely $w_1\cdots w_N
\left(x^{2k}\right)^\circ_L$, where $\left(x^{2k}\right)^\circ_L$
commutes with $x$.  Since $(w_1\cdots w_N)^{-1} y (w_1\cdots w_N) =
x$, it then follows that $\mathbf c_{m-M}(y) = x$. Since $x$ is left
rigid, and $m-M$ is a multiple of $2r$, we finally obtain
$\Psi(y)=x$, that is, $\Psi(\Phi(x))=x$, as we wanted to show.
\end{proof}

In order to finish the proof of Theorem~\ref{T:rigid_USG}, it just
remains to show that the map $\Psi$ can be extended to the arrows of
$USG_R(x)$, so that $\Psi\circ \Phi =\mbox{id}_{USG_L(x)}$. We will
use the following result:

\begin{lemma}\label{L:conjugating_elements_Phi_Psi}
Let $x\in B_n$ be a left rigid braid with $\ell(x)=r>1$. Let $T=2rt$
be such that $\Phi(x)=\mathbf c^{T}_R(x)$ and $\Psi(\Phi(x))=\mathbf
c^T_L(\Phi(x))$. Let $C_T',\ldots, C_1'$ be the conjugating elements for
the iterated right cyclings of $x$,  and let  $C_1,\ldots, C_T$ be the
conjugating elements for the iterated left cyclings of $\Phi(x)$.
That is,
$$
   \Phi(x)= (C_1'\cdots C_T') \: x \: (C_1'\cdots C_T')^{-1}
$$
and
$$
   \Psi(\Phi(x)) = (C_1\cdots C_T)^{-1} \:\Phi(x)\: (C_1\cdots C_T).
$$
Then $C_1\cdots C_T = C_1'\cdots C_T'$.
\end{lemma}

\begin{proof}
Using the notation in the proof of
Proposition~\ref{P:Phi_bijective_on_vertices}, we notice that the
right normal form of $C_1'\cdots C_T'$ is $(y^{2k})^{\circ}_L
(z_N'\cdots z_1')$ for some $k$, and the left normal
form of $C_1\cdots C_T $ is $(w_1\cdots w_N) (x^{2k})^\circ_L$, where
$k$ is the same as above since the supremum of both elements is
precisely $T$. But notice that $(y^{2k})^{\circ}_L (z_N'\cdots z_1')
= (z_N'\cdots z_1') (x^{2k})^\circ_L = (w_1\cdots w_N)
(x^{2k})^\circ_L$, hence the result follows.
\end{proof}

\begin{proof}[Proof of Theorem~\ref{T:rigid_USG}.]
We define $\Psi:\: USG_R(x)^{op} \rightarrow USG_L(x)$ in the natural
way. For every element $u\in USS_R(x)$, we define $\Psi(u)$ as
above, in the same way as $\Phi$ but using right normal forms, that
is, $\Psi(u)= \overleftarrow{\Phi(\overleftarrow u)}$.   In the case
of the arrows of $USG_R(x)^{op}$, we proceed exactly the same way. If $s$
is a simple element such that $s u s^{-1}=v$ with $u,v\in USS_R(x)$,
that is, if $s$ is an arrow in $USG_R(x)^{op}$ going from $v$ to $u$, we
define $\Psi(s)= \overleftarrow{\Phi(\overleftarrow s)}$, where
$\overleftarrow s$ corresponds to an arrow in $USS_L(\overleftarrow
x)$ going from $\overleftarrow u$ to $\overleftarrow v$.

Let us show that, if $s$ is an arrow in $USG_L(x)$ going from $x$ to
$y$, then $\Psi(\Phi(s))=s$.  First, by construction $\Psi(\Phi(s))$
is a simple element conjugating $\Psi(\Phi(x))=x$ to
$\Psi(\Phi(y))=y$, hence $\Psi(\Phi(s))$ is an arrow in $USG_L(x)$
going from $x$ to $y$. We just need to show that $s$ and
$\Psi(\Phi(s))$ are the same as simple elements.

Let $N=2rt$ be big enough, so that $\Phi(x)=\mathbf c^{N}_R(x)$,
$\Phi(y)=\mathbf c^N_R(y)$, $\Psi(\Phi(x))=\mathbf c^N_L(\Phi(x))$
and $\Psi(\Phi(y))=\mathbf c^N_L(\Phi(y))$.  By
Lemma~\ref{L:conjugating_elements_Phi_Psi}, the product of
conjugating elements (on the left) to go from $x$ to $\Phi(x)$ is
the same as the product of conjugating elements (on the right) to go
from $\Phi(x)$ to $\Psi(\Phi(x))=x$. Denote this product by
$\alpha$.  The same happens with $y$ and $\Phi(y)$, and we denote
the corresponding product by $\beta$.  Hence,
$\Psi(\Phi(s))=\Psi(s^{(N)}_R)= \Psi(\alpha s \beta^{-1}) =
\alpha^{-1}( \alpha s \beta^{-1}) \beta = s$, so the result follows.
\end{proof}

% --
%We remark that, since Gebhardt's left transport preserves left
We remark that, since the left transport preserves left
% -- VG
gcd's, $\Phi$ sends minimal arrows of $USG_L(x)$ to minimal arrows
of $USG_R(x)$. By symmetry, $\Psi$ sends minimal arrows in
$USG_R(x)$ to minimal arrows of $USG_L(x)$. Therefore, we have:

\begin{corollary}
\label{C:minimal_USS_graphs}
Let $x\in B_n$ be a left rigid braid with $\ell(x)>1$. The
restriction of $\Phi$ to $minUSG_L(x)$ is an isomorphism of directed
graphs:  $\Phi: \: minUSG_L(x) \rightarrow minUSG_R(x)^{op}$.
\end{corollary}

% --
\subsubsection{$\Phi$ respects the structure of ultra summit graphs}

It was shown in \cite{BGG2} that the arrows of $minUSG_L(x)$, and
similarly those of $minUSG_R(x)^{op}$, can be partitioned naturally
into two categories, namely \textsl{partial cycling} and
\textsl{partial twisted decycling} components.  In this subsection we
show that the isomorphism $\Phi$ is natural in the sense that it
preserves this decomposition of ultra summit graphs.

\begin{proposition}
\label{P:black_or_grey_L}\textup{\cite{BGG2}}
Let $x\in B_n$ with $\ell(x)>0$ and let $s$ be an arrow in
$minUSG_L(x)$ going from $x$ to $x^s$.  Then at least one of the
following conditions holds:
\begin{enumerate}
\vspace{-\topsep}
\item $s\preccurlyeq \iota_L(x)$
\item $s\preccurlyeq \iota_L(x^{-1})$
\end{enumerate}
\end{proposition}
Notice that $\iota_L(x^{-1}) = \partial(\varphi_L(x))$.

\begin{definition}
\label{D:black_or_grey_L}\textup{\cite{BGG2}}
Let $x\in B_n$ with $\ell(x)>0$ and let $s$ be an arrow in
$USG_L(x)$ going from $x$ to $x^s$.
We call $s$ a {\bf partial left cycling} of $x$ and say that the arrow
$s$ is {\bf black} if $s\preccurlyeq \iota_L(x)$.
We call $s$ a {\bf partial twisted left decycling} of $x$ and say that
the arrow $s$ is {\bf grey} if
$s\preccurlyeq \iota_L(x^{-1}) = \partial(\varphi_L(x))$.
\end{definition}

By symmetry we have
\begin{proposition}
\label{P:black_or_grey_R}\textup{\cite{BGG2}}
Let $x\in B_n$ with $\ell(x)>0$ and let $s$ be an arrow in
$minUSG_R(x)$ going from $x$ to $x^s$.  Then at least one of the
following conditions holds:
\begin{enumerate}
\vspace{-\topsep}
\item $\iota_R(x)\succcurlyeq s$
\item $\iota_R(x^{-1})\succcurlyeq s$
\end{enumerate}
\end{proposition}
Notice that $\iota_R(x^{-1}) = \partial^{-1}(\varphi_R(x))$.

\begin{definition}
\label{D::black_or_grey_R}\textup{\cite{BGG2}}
Let $x\in B_n$ with $\ell(x)>0$ and let $s$ be an arrow in
$USG_R(x)$ going from $x$ to $x^s$.
We call $s$ a {\bf partial right cycling} of $x$ and say that the
arrow $s$ is {\bf black} if $\iota_R(x)\succcurlyeq s$.
We call $s$ a {\bf partial twisted right decycling} of $x$ and say
that the arrow $s$ is {\bf grey} if
$\partial^{-1}(\varphi_R(x)) = \iota_R(x^{-1})\succcurlyeq s$.
\end{definition}

Note that the intuitive meaning of ``cycling'' (respectively
``decycling'') is to move the first simple factor to the end
(respectively, the last simple factor to the front) with respect to
the normal form under consideration.  Note also that $\tau\circ \mathbf d_L(x)
= \tau(x^{\varphi_L(x)^{-1}})=\tau(x^{\iota_L(x^{-1})\Delta^{-1}})
= x^{\iota_L(x^{-1})}$ and that
$\tau^{-1}\circ \mathbf d_R(x)
= \tau^{-1}(x^{\varphi_R(x)})=\tau^{-1}(x^{\iota_R(x^{-1})^{-1}\Delta})
= x^{\iota_R(x^{-1})^{-1}}$
Hence, the definitions of ``partial cycling'' and ``partial twisted
decycling'' are natural:  a partial cycling or decycling corresponds
to moving a prefix or suffix of the first or last simple factor;
``twisting'' refers to composition with $\tau$.

Partial cyclings and partial twisted decyclings are preserved by
the graph isomorphism $\Phi$ according to the following results.

\begin{proposition}
\label{P:Phi_preserves_colours}
Let $x\in B_n$ be a rigid braid with $\ell(x)>1$, and let $s$ be an arrow from $x$ to $y$ in $USG_L(x)$
such that $s\preccurlyeq\iota_L(x)$.  Then, $\Phi(s)$ is an arrow from
$\Phi(y)$ to $\Phi(x)$ in $USG_R(x)$ such that
$\iota_R(\Phi(y))\succcurlyeq \Phi(s)$.
\end{proposition}
\begin{proof}
Recall that $\Phi(s)$ is defined via iterated transport.  As transport
is monotonic, we obtain $\Phi(s)\preccurlyeq\iota_L(\Phi(x))$.
Moreover, $\Phi(s)$ is simple.
If $\Phi(x)=\Delta^p x_1\cdots x_r$ is in left normal form, we have
$\tau^p(\Phi(s))\preccurlyeq x_1$, whence
$\Phi(y) = \Phi(s)^{-1}\Phi(x)\Phi(s)
  = \Delta^p(\tau^p(\Phi(s))^{-1}x_1)x_2\cdots x_r\Phi(s)$.
The latter implies $\iota_R(\Phi(y))\succcurlyeq\Phi(s)$ as claimed,
since $\inf(\Phi(y))=\inf(\Phi(x))=p$.
\end{proof}

\begin{corollary}
$\Phi$ and $\Psi$ are isomorphisms of directed graphs preserving the colours of arrows.
\end{corollary}
\begin{proof}
We know that $\Phi$ and $\Psi$ are isomorphisms of directed graphs by
Theorem~\ref{T:rigid_USG};  it remains to be shown that they
preserve the colours of arrows.

By Proposition~\ref{P:Phi_preserves_colours}, the image of a black
arrow under $\Phi$ is a black arrow. Applying Proposition~\ref{P:Phi_preserves_colours} to $x^{-1}$, which is also a rigid element with $\ell(x^{-1})>1$, it follows that the image of a grey arrow under $\Phi$ is a grey arrow. The analogous result holds for $\Psi$ by symmetry. 
\end{proof}
% -- VG

\vspace{.3cm}
\noindent {\footnotesize
\begin{minipage}[t]{5.4cm}
{\bf Juan Gonz\'alez-Meneses:} \\
Dept. \'{A}lgebra.  Facultad de Matem\'{a}ticas\\
Universidad de Sevilla. \\
Apdo. 1160. \\
41080 Sevilla (SPAIN)
\\ E-mail:  meneses@us.es
\\ URL: www.personal.us.es/meneses
\end{minipage}
\hfill
\begin{minipage}[t]{5.2cm}
{\bf Volker Gebhardt:} \\
School of Computing and Mathematics\\
University of Western Sydney \\
Locked Bag 1797 \\
Penrith South DC NSW 1797, Australia
\\ E-mail: v.gebhardt@uws.edu.au
\end{minipage}
}

\end{document}